\documentclass[preprint,3p,times]{elsarticle}

\usepackage{amsmath}
\usepackage{amssymb}
\usepackage{amsthm}
\usepackage{graphicx}
\usepackage{natbib}

\usepackage{moreverb}

\usepackage{xfrac} % for \sfrac{1}{2}
\usepackage{siunitx}

\usepackage{algorithm}
\usepackage[noend]{algpseudocode}
\newtheorem{theorem}{Theorem}[section]

\usepackage{booktabs}
\usepackage{xcolor}

\usepackage[capitalise,noabbrev]{cleveref}
\crefformat{appendix}{#2#1#3}
\crefformat{equation}{(#2#1#3)}

\renewcommand{\vec}[1]{\boldsymbol{#1}}

\newcommand{\ssff}{FF}
\newcommand{\ssfc}{FC}
\newcommand{\sscf}{CF}
\newcommand{\sscc}{CC}
\newcommand{\ssf}{F}
\newcommand{\ssc}{C}

\newcommand{\ndof}{n_{\text{per DoF}}}
\newcommand{\ncycle}{n_{\text{per DoF per cycle}}}

\graphicspath{{./figures/}}

\journal{Journal of Computational and Applied Mathematics}

\begin{document}

\begin{frontmatter}

\newcommand{\thistitle}{Coarse-Grid Selection Using Simulated Annealing}
\title{\thistitle}

\author[scomp-mun]{T. U. Zaman\corref{cor}}
\cortext[cor]{Corresponding author.}
\ead{tzaman@mun.ca}

\author[math-mun]{S. P. MacLachlan}
%\author[math-mun]{S. P. MacLachlan\fnref{fn-spm}}
%\fntext[fn-spm]{Present address: Department of Mathematics and Statistics,
%Memorial University of Newfoundland,
%NL, Canada.}
%\ead{smaclachlan@mun.ca}

\author[comp-uiuc]{L. N. Olson}
%\author[comp-uiuc]{L. N. Olson\fnref{fn-lo}}
%\fntext[fn-lo]{Present address: Department of Computer Science,
%University of Illinois at Urbana-Champaign, Illinois, USA.}
%\ead{lukeo@illinois.edu}

\author[mech-uiuc]{M. West}
%\author[mech-uiuc]{M. West\fnref{fn-mw}}
%\fntext[fn-mw]{Present address: Department of Mechanical Science and Engineering,
%University of Illinois at Urbana-Champaign, Illinois, USA.}
%\ead{mwest@illinois.edu}

\address[scomp-mun]{Scientific Computing Program,
Memorial University of Newfoundland, NL, Canada.}
\address[math-mun]{Department of Mathematics and Statistics,
Memorial University of Newfoundland, NL, Canada.}
\address[comp-uiuc]{Department of Computer Science,
University of Illinois at Urbana-Champaign, Illinois, USA.}
\address[mech-uiuc]{Department of Mechanical Science and Engineering,
University of Illinois at Urbana-Champaign, Illinois, USA.}

\begin{abstract}
Multilevel techniques are efficient approaches for solving the large linear systems that arise from discretized partial differential equations and other problems. While geometric multigrid requires detailed knowledge about the underlying problem and its discretization, algebraic multigrid aims to be less intrusive, requiring less knowledge about the origin of the linear system. A key step in algebraic multigrid is the choice of the coarse/fine partitioning, aiming to balance the convergence of the iteration with its cost. In work by MacLachlan and Saad~\cite{maclachlan2007greedy}, a constrained combinatorial optimization problem is used to define the ``best'' coarse grid within the setting of a two-level reduction-based algebraic multigrid method and is shown to be NP-complete. Here, we develop a new coarsening algorithm based on simulated annealing to approximate solutions to this problem, which yields improved results over the greedy algorithm developed previously. We present numerical results for test problems on both structured and unstructured meshes, demonstrating the ability to exploit knowledge about the underlying grid structure if it is available.
\end{abstract}


%\jnlcitation{\cname{%
    %\author{T. U. Zaman},
    %\author{S. P. MacLachlan},
    %\author{L. N. Olson}, and
    %\author{M. West}} (\cyear{2021}),
  %\ctitle{\thistitle}, \cjournal{Numerical Linear Algebra with Applications}, \cvol{Volume 1}.}

%\maketitle

%\footnotetext{\textbf{Abbreviations:} AMG, algebraic multigrid; AMGr, reduction-based AMG\@; SA, simulated annealing}

\end{frontmatter}

%%%%%%%%%%%%%%%%%%%%%%%%%%%%%%%%%%%%%%%%%%%%%%%%%%%%%%%%%%%%%%%%%%%%%%%%%%%%%%

\section{Introduction}

Multigrid and other multilevel methods are well-established
algorithms for the solution of large linear systems of
equations that arise in many areas of computational science and
engineering.  Multigrid methods arise from the
observation that basic iterative methods, such as the (weighted)
Jacobi and Gauss-Seidel iterations, effectively eliminate only part of
the error in an approximation to the solution of the problem, and that
the complementary space of errors can be corrected from a coarse
representation of the problem.  While geometric multigrid has
been shown to be highly effective for problems where a hierarchy of
discrete models can be built directly, many problems benefit from the
use of algebraic multigrid (AMG) techniques, where graph-based
algorithms and other heuristics are used to define the multigrid
hierarchy directly from the matrix of the finest-grid problem.

First introduced in the early 1980's~\cite{ABrandt_SFMcCormick_JWRuge_1982,
  ABrandt_SFMcCormick_JWRuge_1984a}, AMG is now widely used and
  available in standard libraries, such as hypre~\cite{RDFalgout_UMYang_2002a, VEHenson_UMYang_2002a}, PETSc~\cite{SBalay_etal_2008a, SBalay_etal_2010a}, Trilinos~\cite{MAHeroux_etal_2005a, MWGee_etal_2006a}, and PyAMG~\cite{pyamg}.  While there are many
variants of AMG (as discussed below), the common
features of AMG algorithms are the use of graph-based algorithms to
construct a coarse grid from the given fine-grid discretization matrix
(possibly with some additional information, such as geometric
locations of the degrees of freedom for elasticity problems~\cite{vanvek1996algebraic}) followed by construction of an algebraic
interpolation operator from the coarse grid to the fine grid.  Both of
these components have received significant attention in the
literature, with an abundance of schemes for creating the coarse grid
and for determining the entries in interpolation.  In this paper, we
focus on the coarse-grid correction process, adopting a combinatorial
optimization viewpoint on the selection of coarse-grid points in the
classical AMG setting.

One of the primary differences between so-called classical
(Ruge-St\"uben) AMG and (smoothed) aggregation approaches is in how
the coarse grid itself is defined~\cite{KStuben_2001a}.
In aggregation-based approaches, coarse grids are defined by clustering
fine-grid points into aggregates (or subdomains), while, in classical AMG,
a subset of the fine-grid degrees of freedom (points) is selected in order to form
the coarse grid.  In the original method~\cite{ABrandt_SFMcCormick_JWRuge_1984a, JWRuge_KStuben_1987a}, this
was accomplished using a greedy algorithm to construct a maximal
independent set over the fine-grid degrees of freedom as a ``first
pass'' at forming the coarse-grid set, which is then augmented by a
``second pass'' algorithm that ensures suitable interpolation can be
defined from the final coarse grid.  Many alternative strategies to
forming the coarse grid have been proposed over the years,
particularly for the parallel case, to overcome the inherent
sequentiality of the original algorithm~\cite{DMAlber_LNOlson_2007a,
  AJCleary_RDFalgout_VEHenson_JEJones_1998a, VEHenson_UMYang_2002a,
  HDeSterck_UMYang_JJHeys_2006a, HDeSterck_etal_2008b}.
Other approaches, for example based on redefining the notion of
strength of connection~\cite{LNOlson_etal_2010a, JBrannick_etal_2006a,
EChow_2003a, OBroeker_2003a} or use of compatible relaxation
principles~\cite{OLivne_2004a, JJBrannick_RDFalgout_2010a}, have also
been proposed.

Much of this work has been motivated by the failure of the classical
AMG heuristics to yield robust solvers for difficult classes of
problems, such as the Helmholtz equation, convection-dominated flows,
or coupled systems of PDEs.  While much research has tried to
generalize these heuristics to give sensible algorithms for broader classes
of problems, another approach is to consider methods that abandon
these heuristics in favour of algorithms that are directly tied to
rigorous multigrid convergence theory.  In recent years, two main classes
of AMG algorithms have arisen in this direction.  One class of methods
arises in the aggregation setting, where well-chosen pairwise
aggregation algorithms coupled with suitable relaxation can yield
convergence guarantees for symmetric, diagonally dominant
M-matrices~\cite{doi:10.1137/100818509, doi:10.1137/120876083}.  The
second class of methods are those based on
reduction-based AMG principles, where there is a direct connection
between properties of the fine-grid submatrix and the guaranteed
convergence rate of the scheme~\cite{SPMacLachlan_TAManteuffel_SFMcCormick_2006a,
  maclachlan2007greedy, SMacLachlan_YSaad_2007b,
  FGossler_RNabben_2016a}.  Reduction-based multigrid methods~\cite{MRies_UTrottenberg_GWinter_1983a} are a generalization of cyclic
reduction~\cite{PNSwarztrauber_1977a}, in which conditions on the coarse/fine
partitioning are relaxed so that while the exact Schur
complement may not be sparse, it can be accurately approximated by a
sparse matrix.  This notion is made rigorous by MacLachlan et al.~\cite{SPMacLachlan_TAManteuffel_SFMcCormick_2006a}, who propose a
theoretical framework to analyze convergence based on diagonal
dominance of the fine-grid submatrix (in contrast to the fine-grid
submatrix being diagonal in the classical setting of cyclic
reduction), and extended to Chebyshev-style relaxation schemes by Gossler
and Nabben~\cite{FGossler_RNabben_2016a}.  While the theoretical
guarantees offered by such methods are attractive, these approaches
suffer from two key limitations.  First, the classes of problems
for which convergence rates are guaranteed are limited
(requiring diagonal dominance and/or M-matrix structure).  Second,
the guaranteed convergence rates of stationary iterations are
generally poor in comparison to methods based on classical
heuristics.  An active area of research (disjoint from the goals of
this paper) is addressing the first
limitation~\cite{doi:10.1137/18M1193761}, while the second can be
ameliorated by the use of Krylov subspace methods to accelerate
stationary convergence.  We note that there are
other families of theoretical analysis of AMG
algorithms~\cite{SMacLachlan_LOlson_2014a}, including theory tailored
to heuristic approaches, such as compatible
relaxation~\cite{JJBrannick_RDFalgout_2010a}.  Much of
this theory, however, is of a confirmational and not a predictive nature, i.e.,
the convergence bounds rely on properties of the output of heuristics
for choosing coarse grids and interpolation operators, and not on the
inputs to these processes, such as the system matrix and algorithmic
parameters.  Thus, while it provides guidance in the development of
heuristic methods, it does not provide the concrete convergence bounds
offered by pairwise-aggregation or AMGr methods.

While the pairwise
aggregation methodology~\cite{doi:10.1137/100818509, doi:10.1137/120876083} provides
practical algorithms to generate coarse grids that satisfy their
convergence guarantees, this is a notable omission
in the work of MacLachlan et al.~\cite{SPMacLachlan_TAManteuffel_SFMcCormick_2006a} and much of
the work on AMGr.
MacLachlan and Saad~\cite{maclachlan2007greedy, SMacLachlan_YSaad_2007b}
identify that the choice of optimal coarse grids for AMGr
can be quantified as the solution of an NP-complete integer linear
programming problem.  They use this formulation to
guide construction of a family of greedy coarsening algorithms that aim to
maximize the size of the fine-grid set subject to maintaining a fixed
level of diagonal dominance in the fine-grid submatrix (and,
consequently, the resulting convergence bound on the AMG method).
While the greedy algorithm was tested on a range of
problems in~\cite{maclachlan2007greedy}, these problems were limited
primarily to bilinear finite-element
discretization on structured grids, with only a few matrices outside this class.  When we assessed its performance on a broader class of
problems, we exposed some simple test cases where it dramatically
fails to perform well, such as standard five-point finite difference
discretization, motivating further work in this area.
In this paper, we follow the theory and framework for AMGr~\cite{SPMacLachlan_TAManteuffel_SFMcCormick_2006a, maclachlan2007greedy}, but aim to overcome some limitations of the underlying {greedy} coarsening
algorithm. In particular, we develop a new coarsening
algorithm based on simulated annealing to partition the coarse and fine points.
Numerical results show that this algorithm achieves smaller coarse
grids (than the greedy approach)
that satisfy
the same diagonal-dominance criterion.  Hence, these grids lead to
moderately more efficient (two-level) algorithms, by reducing the size of the
coarse-grid problem without changing the convergence bound guaranteed
by the theory.  We emphasize that this work is presented more as a
proof-of-concept than as a practical coarsening algorithm, due to the
high cost of simulated annealing.  In related work, we have shown that
the ``online'' cost of the approach proposed here can be traded for a
significant ``offline'' cost using a reinforcement learning
approach~\cite{ATag_etal_2021a}.  However, a necessary next step in the work
proposed here is in investigating more cost-effective alternatives to
these methods.

This paper is arranged as follows. In~\cref{sec:amgcoarsening-cgsa}, AMG coarsening
and the existing greedy
coarsening algorithm for AMGr are discussed. The new coarsening algorithm based on simulated
annealing is outlined in~\cref{sec:sacoarsening-cgsa}. Numerical experiments are given in~\cref{sec:results-cgsa},
demonstrating performance of this approach on both
isotropic and anisotropic problems, discretized on both structured and
unstructured meshes.
Concluding remarks and potential future research are discussed
in~\cref{sec:conclusion-cgsa}.

%%%%%%%%%%%%%%%%%%%%%%%%%%%%%%%%%%%%%%%%%%%%%%%%%%%%%%%%%%%%%%%%%%%%%%%%%%%%%%
\section{AMG coarsening}\label{sec:amgcoarsening-cgsa}

Geometric multigrid is known to be highly effective
for many problems discretized on structured meshes. However, it is
naturally more difficult
to make effective coarse grids (and effective coarse-grid operators) for problems discretized on unstructured meshes. AMG was
developed specifically to address this, automating the formation of
coarse-grid matrices without any direct mesh information.
The first coarsening approach in AMG, often called classical or Ruge-St\"{u}ben (RS)
coarsening, was developed by Brandt, McCormick, and Ruge~\cite{ABrandt_SFMcCormick_JWRuge_1982, ABrandt_SFMcCormick_JWRuge_1984a} and
later by Ruge and St\"{u}ben~\cite{JWRuge_KStuben_1987a}.  We
summarize this approach below, primarily to allow comparison with the
reduction-based AMG method~\cite{SPMacLachlan_TAManteuffel_SFMcCormick_2006a,
  MRies_UTrottenberg_GWinter_1983a} and the greedy coarsening approach proposed by MacLachlan and Saad~\cite{maclachlan2007greedy} that is the starting point for the
research reported herein.

AMG algorithms make use of a two-stage approach, where
a setup phase  precedes the cycling in the solve
phase.  Common steps in AMG setup phase include recursively choosing a coarse
grid and defining intergrid transfer
operators.
The details of these processes depend on the specifics of
the AMG algorithm under consideration, with both point-based and
aggregation-based approaches to determining the coarse grid, and a
variety of interpolation schemes possible for both of these.

%%%%%%%%%%%%%%%%%%%%%%%%%%%%%%%%%%%%%%%%%%%%%%%%%%%%%%%%%%%%%%%%%%%%%%%%%%%%%%
\subsection{Classical coarsening}

As in geometric multigrid,
the coarse grid in AMG is selected so that errors
not reduced by relaxation
can be accurately approximated on coarse grids. In an effective
scheme, these errors are
interpolated accurately from coarse grids that have
substantially fewer degrees of freedom than the next finer grid, thus significantly reducing
the cost of solving the coarse-grid residual problem.

In classical AMG, for an $n\times n$ matrix $A$, the index set $\Omega =
\{1,\dots,n\}$ is split into sets $C$ and $F$, with $\Omega = C \cup F$ and $C
\cap F =\emptyset$.  Each degree of freedom (DoF) or point $i\in C$ is denoted a $C$-point and a point
$j\in F$ is denoted an $F$-point.
This splitting is constructed by considering the so-called \textit{strong} connections in the graph of matrix between the fine-level
variables.
Given a threshold value,
$0<\theta \leq 1$, the variable $u_{j}$ is said to
strongly influence $u_{i}$ if
$-A_{ij} \geq \theta   \displaystyle{\max_{k \neq i}}\{-A_{ik}\}$, where
$A_{ij}$ is the coefficient of $u_{j}$ in the $i$th equation. The set of
points that strongly influence $i$, denoted by $S_{i}$, is defined as the
set of points on which point $i$ strongly depends. The set of points that
are strongly influenced by $i$ is denoted by $S_{i}^T$. Two heuristics are followed
in RS coarsening to select a coarse grid:
\begin{description}
\item[\textit{H1}]:  For every $F$-point, $i$, every point $j\in
  S_{i}$ should either be a coarse-grid point or should strongly
  depend on at least one point in $C$ that also strongly influences $i$.
  \item[\textit{H2}]: The set
of $C$-points should be a maximal subset of all points, where no $C$-point strongly depends on another $C$-point.
\end{description}
%\textit{H1} is intended to ensure that the fine-grid variables are able to be accurately
%interpolated
%from the coarse-grid variables using classical AMG
%interpolation. On the other hand, \textit{H2} constrains the size of set $C$,
%in effect controlling computational cost.
In practice, strong enforcement of both of \textit{H1} and \textit{H2}
is not always possible; the classical interpolation formula relies on
\textit{H1} being strongly enforced, while \textit{H2} is used only to
encourage the selection of small, sparse coarse grids.

%%%%%%%%%%%%%%%%%%%%%%%%%%%%%%%%%%%%%%%%%%%%%%%%%%%%%%%%%%%%%%%%%%%%%%%%%%%%%%
\subsection{Greedy coarsening and underlying optimization}

While the classical coarsening algorithm has proven effective for many
problems when coupled with suitable construction of the interpolation
operator~\cite{JWRuge_KStuben_1987a}, the process provides few guarantees in practice.
Indeed, there is an inherent disconnect
between the selection of any single coarse-grid point
and the impact on the quality of the resulting interpolation operator.
This has motivated coupled approaches to selecting coarse grids and defining
interpolation, including
\textit{reduction-based} AMG, or AMGr.

Reduction-based multigrid was proposed by Ries et
al.~\cite{MRies_UTrottenberg_GWinter_1983a}, building on earlier work
aiming to improve multigrid convergence for the standard
finite-difference (FD) Poisson problem.  The fundamental idea of
reduction-based multigrid lies in defining the multigrid components to
approximate those of cyclic-reduction
algorithms~\cite{PNSwarztrauber_1977a}.  In particular, one way to
interpret cyclic reduction is as a two-level multigrid method with
idealized relaxation, interpolation, and restriction operators.
Suppose that the coarse/fine partitioning is already determined, and
consider the reordering of the linear system $A\vec{x} = \vec{b}$ according to the partition, writing
\[
A=
\begin{bmatrix}
 A_{\ssff} & -A_{\ssfc} \\
-A_{\sscf} &  A_{\sscc}
\end{bmatrix}
\quad
\vec{x} =
\begin{bmatrix}
  \vec{x}_{\ssf}\\
  \vec{x}_{\ssc}
\end{bmatrix}
\quad
\vec{b} =
\begin{bmatrix}
  \vec{b}_{\ssf}\\
  \vec{b}_{\ssc}
\end{bmatrix}.
\]
An \textit{exact}
algorithm for the solution of $A\vec{x} = \vec{b}$ in this partitioned
form is given by
\begin{enumerate}
  \item $\vec{y}_{\ssf} = A_{\ssff}^{-1}\vec{b}_\ssf$,
  \item Solve $\left(A_{\sscc} - A_{\sscf}A_{\ssff}^{-1}A_{\ssfc}\right)\vec{x}_{\ssc}
    = \vec{b}_{\ssc} + A_{\sscf}\vec{y}_{\ssf}$,
  \item $\vec{x}_{\ssf} = \vec{y}_{\ssf} + A_{\ssff}^{-1}A_{\ssfc}\vec{x}_{\ssc}$.
\end{enumerate}
This can be turned into an iterative method for solving
$A\vec{x}=\vec{b}$ in the usual way, by replacing the right-hand side
vector, $\vec{b}$, by the evolving residual, and computing updates to a
current approximation, $\vec{x}^{(k)}$, giving
\begin{enumerate}
  \item $\vec{x}^{(k+\sfrac{1}{2})}_{\ssf} = \vec{x}^{(k)}_{\ssf} +
    A_{\ssff}^{-1}\left(\vec{b}_{\ssf} - A_{\ssff}\vec{x}^{(k)}_{\ssf} + A_{\ssfc}\vec{x}^{(k)}_{\ssc}\right)$,
  \item Solve $\left(A_{\sscc} - A_{\sscf}A_{\ssff}^{-1}A_{\ssfc}\right)\vec{y}_{\ssc}
    = \vec{b}_{\ssc} + A_{\sscf}\vec{x}^{(k+\sfrac{1}{2})}_{\ssf} - A_{\sscc}\vec{x}^{(k)}_{\ssc}$,
  \item $\vec{x}^{(k+1)}_{\ssc} = \vec{x}^{(k)}_{\ssc} + \vec{y}_{\ssc}$,
  \item $\vec{x}^{(k+1)}_{\ssf} = \vec{x}^{(k+\sfrac{1}{2})}_{\ssf} + A_{\ssff}^{-1}A_{\ssfc}\vec{y}_{\ssc}$.
\end{enumerate}
In this form, this remains an exact algorithm: given any
initial guess, $\vec{x}^{(0)}$, we have the solution $\vec{x} =
\vec{x}^{(1)}$.  A truly iterative method results from approximating
the three instances of $A_{\ssff}^{-1}$ in the above algorithm, namely
\begin{equation}
  \hat{A}_{\ssff}^{-1} \approx A_{\ssff}^{-1}
  \qquad
  \hat{A}_{\ssc} \approx A_{\sscc} - A_{\sscf}A_{\ssff}^{-1}A_{\ssfc}
  \qquad
  W_{\ssfc} \approx A_{\ssff}^{-1}A_{\ssfc},
\end{equation}
leading to the following iteration:
\begin{enumerate}
  \item $\vec{x}^{(k+\sfrac{1}{2})}_{\ssf} = \vec{x}^{(k)}_{\ssf} +
    \hat{A}_{\ssff}^{-1}\left(\vec{b}_{\ssf} - A_{\ssff}\vec{x}^{(k)}_{\ssf} + A_{\ssfc}\vec{x}^{(k)}_{\ssc}\right)$,
  \item Solve $\hat{A}_{\ssc}\vec{y}_{\ssc}
    = \vec{b}_{\ssc} + A_{\sscf}\vec{x}^{(k+\sfrac{1}{2})}_{\ssf} - A_{\sscc}\vec{x}^{(k)}_{\ssc}$,
  \item $\vec{x}^{(k+1)}_{\ssc} = \vec{x}^{(k)}_{\ssc} + \vec{y}_{\ssc}$,
  \item $\vec{x}^{(k+1)}_{\ssf} = \vec{x}^{(k+\sfrac{1}{2})}_{\ssf} + W_{\ssfc}\vec{y}_{\ssc}$.
\end{enumerate}
In multigrid terminology, the first step can be seen as a so-called
$F$-relaxation step, where we approximate $A_{\ssff}^{-1}$ using some
simple approximation, such as with weighted Jacobi or Gauss-Seidel.
The second step is a coarse-grid correction step, where the residual
after $F$-relaxation is restricted by injection, and the coarse-grid
correction, $\vec{y}_{\ssc}$, is computed by solving a linear system with
an approximation, $\hat{A}_{\ssc}$, of the true Schur complement, $A_{\sscc} -
A_{\sscf}A_{\ssff}^{-1}A_{\ssfc}$.  The final two steps correspond to an
interpolation of the coarse-grid correction, writing the interpolation
operator $P = \left[\begin{smallmatrix} W_{\ssfc}
    \\ I \end{smallmatrix}\right]$, for some approximation of the
    ideal interpolation operator, $W_{\ssfc} \approx A_{\ssff}^{-1}A_{\ssfc}$.  If
the approximation in the first step is sufficiently poor that a
significant residual is expected to remain at the $F$-points after
relaxation, it is also reasonable to augment the restriction in the
second step, using either the transpose of interpolation or some
better approximation to an ideal restriction operator than just injection.

As written above, the algebraic form of reduction-based multigrid retains the key disadvantage
of classical AMG~---~there is little connection between the algorithmic
choices (of the coarse/fine partitioning, and the approximations of
$A_{\ssff}$, of the Schur complement, and of the idealized interpolation and
restriction operators) with the resulting convergence of the
algorithm.
To address these limitations, MacLachlan et al.~\cite{SPMacLachlan_TAManteuffel_SFMcCormick_2006a} proposed a reduction-based
AMG algorithm (AMGr) that attempts to connect convergence with
properties of $A_{\ssff}$.
In particular, they presented a two-level
convergence theory with a convergence rate quantified by the
approximation of $A_{\ssff} \approx D_{\ssff}$, with the expectation that
application of $D_{\ssff}^{-1}$ to both a vector and to $A_{\ssfc}$ is
computationally feasible.  In what follows, we use the standard
notation that matrices $A \succeq B$ when $\vec{x}^T A\vec{x} \succeq
\vec{x}^T B\vec{x}$ for all vectors $\vec{x}$.

\begin{theorem}\label{thm:amgr-cgsa}\cite{SPMacLachlan_TAManteuffel_SFMcCormick_2006a}
Consider the symmetric and positive-definite matrix
$A=\big[ \begin{smallmatrix}
    A_{\ssff} & -A_{\ssfc} \\
    -A_{\ssfc}^T & A_{\sscc}
\end{smallmatrix}\big]$ such that $A_{\ssff}=D_{\ssff}+\mathcal{E}$, with $D_{\ssff}$ symmetric,
$0\preceq \mathcal{E} \preceq \epsilon D_{\ssff}$, and
$\big[ \begin{smallmatrix}
    D_{\ssff} & -A_{\ssfc} \\
    -A_{\ssfc}^T & A_{\sscc}
\end{smallmatrix}\big]
\succeq 0$, for some $\epsilon \geq 0$.
Define relaxation with error-propagation operator $R=
\big(I-\sigma \big( \begin{smallmatrix}
      D_{\ssff}^{-1} & 0 \\
0 & 0
\end{smallmatrix}\big)A \big)$
for $\sigma=2/(2+\epsilon)$, interpolation $P=\big[
\begin{smallmatrix}
  D_{\ssff}^{-1}A_{\ssfc} \\
I
\end{smallmatrix}\big]$, and coarse-level correction with
error-propagation operator $T=I-P(P^T A P)^{-1}P^T A$. Then the
multigrid cycle with $\nu$ pre-relaxation sweeps, coarse-level correction,
and $\nu$ post-relaxation sweeps has error propagation operator
$\text{MG}_2=R^{\nu}\cdot T \cdot R^{\nu}$ which satisfies
\begin{equation}\label{eqn:amgr_bound-cgsa}
||\text{MG}_2||_A\leq \left(\frac{\epsilon}{1+\epsilon} \left(1+\left(\frac
{\epsilon^{2\nu -1}}{(2+\epsilon)^{2\nu}} \right)\right)\right)^{1/2} < 1.
\end{equation}
\end{theorem}

If a partitioning and approximation, $D_{\ssff}$, of $A_{\ssff}$ are found
that satisfy the assumptions given above, then this theorem
establishes existence of an interpolation
operator, $P$, giving multigrid convergence with a direct tie to the
approximation parameter,~$\epsilon$. In this theory,
tight spectral equivalence between $D_{\ssff}$ and $A_{\ssff}$ is required
to ensure good performance of the solver. This leads to the goal of
constructing the fine-grid set so that there is a guarantee of tight spectral
equivalence between $D_{\ssff}$ and $A_{\ssff}$.

To meet that goal, MacLachlan and Saad~\cite{maclachlan2007greedy}
proposed to partition the rows and columns of the matrix $A$
in order to ensure the diagonal dominance
of $A_{\ssff}$, so that $D_{\ssff}$ could be chosen as a diagonal matrix. In
particular, a diagonal dominance factor is introduced for
each row $i$, defined as
\[
  \theta_i=\frac{|A_{ii}|}{\sum_{j\in F}|A_{ij}|},
\]
With this, $A_{\ssff}$ is said to be $\theta$-diagonally dominant if $\theta_i \geq \theta$
for all $i \in F$, where $\theta>1/2$ measures the diagonal dominance
of $A_{\ssff}$. If $A_{\ssff}$ is $\theta$-diagonally
dominant, then it is shown that the diagonal matrix, $D_{\ssff}$, with
$(D_{\ssff})_{ii} = (2-\frac{1}{\theta})a_{ii}$ for all $i\in F$ leads
to $0\preceq \mathcal{E} \preceq
\frac{2-2\theta}{2\theta-1}D_{\ssff}$, giving a $\theta$-dependent
convergence bound if the other assumptions of~\Cref{thm:amgr-cgsa} are
satisfied.  Furthermore, if $A$ is symmetric, positive-definite, and
diagonally dominant, then this prescription for $D_{\ssff}$ guarantees
that all conditions of~\Cref{thm:amgr-cgsa} are satisfied, so long as
$A_{\ssff}$ is $\theta$-diagonally dominant.

In
addition to establishing this connection between the diagonal dominance
parameter $\theta$ and the convergence parameter, $\epsilon$,
MacLachlan and Saad~\cite{maclachlan2007greedy} posed the partitioning
algorithm as an
optimization problem for given $\theta>1/2$, asking for the largest F-set such that
$\theta_i\geq \theta$ for every
$i \in F$. Such a $\theta$-dominant $A_{\ssff}$
guarantees good
equivalence between a diagonal matrix, $D_{\ssff}$, and $A_{\ssff}$, and the largest
such $F$-set would make the coarse-grid problem smallest.
This leads to an optimization problem of the form
\begin{equation}\label{eq:linear_program-cgsa}
  \begin{split}
\max_{F\subset\Omega}   & \lvert F\rvert,\\
\text{subject to } & |A_{ii}| \geq \theta \sum_{j \in F} |A_{ij}|, \,  \forall i \in F.
  \end{split}
\end{equation}

The greedy algorithm~\cite{maclachlan2007greedy} for~\cref{eq:linear_program-cgsa}
selects rows and
columns of the matrix $A$ for
$A_{\ssff}$ in a greedy manner to ensure the diagonal dominance
of $A_{\ssff}$, as described in~\cref{alg:greedy_alg-cgsa}. Here, the set $U$ contains all the DoFs
that are unpartitioned (not yet assigned to be coarse or fine points). Initially,
all degrees of freedom are assigned to $U$, while the $F$ and $C$ sets are
empty. The rows that directly satisfy the diagonal dominance criterion
are initially added to the $F$-set (and removed from $U$). If there
are no diagonally dominant DoFs in the $U$ set, then the point
having the least diagonal dominance is made a $C$-point. This
selection may make some other points in
$U$ diagonally dominant, whereupon these DoFs are added to the
$F$-set and removed from $U$, and the process is repeated until the
set $U$ becomes empty.
\begin{algorithm}[!ht]
\caption{Greedy coarsening algorithm}\label{alg:greedy_alg-cgsa}
\begin{algorithmic}[1]
  \Function{greedy-coarsening}{$A$, $\theta$}
  \State $U \gets \{1,2,\ldots,n\}$, $F \gets \emptyset$, $C \gets \emptyset$
  \For {$i \gets 1,\ldots,n$}
      \State $\hat{\theta}_i \gets \frac{\lvert A_{ii}\rvert}{\sum_{j \in F\cup U}\lvert A_{ij}\rvert}$
      \If {$\hat{\theta}_i \geq \theta$}
          \State $F \gets F\cup \{i\}$, $U \gets U\setminus \{i\}$
      \EndIf
  \EndFor
  \While {$U\neq \emptyset$}
      \State $j \gets \text{argmin}_{i\in U}\{\hat{\theta}_i\}$
      \State $U \gets U \setminus \{j\}$, $C \gets C \cup \{j\}$
      \For {$i \in U \cap \text{Adj}(j)$} \Comment{$\text{Adj}(j)=\{k:A_{jk}\neq 0\}$}
          \State $\hat{\theta}_i\gets \frac{\lvert A_{ii}\rvert}{\sum_{k\in F\cup U}\lvert A_{ik}\rvert}$
          \If {$\hat{\theta}_i \geq \theta$}
              \State $F \gets F \cup \{i\}$, $U \gets U \setminus \{i\}$
          \EndIf
      \EndFor
      \EndWhile
      \State Return $F$, $C$
\EndFunction
\end{algorithmic}
\end{algorithm}

We note that there is a slightly counter-intuitive connection between the
quality of solution to \cref{eq:linear_program-cgsa} and the convergence
of the multigrid method.  ``Bad'' solutions to
\cref{eq:linear_program-cgsa}, meaning those with satisfied constraints
but that are far from optimality, have much bigger $C$-sets than ``good'' solutions do and, consequently, the resulting two-level convergence can be much better than that
guaranteed by the theoretical bounds.  Extreme examples of this occur when the set $F$ is an independent set, so that $A_{FF}$ is a diagonal matrix, and the resulting two-grid method is exact.  While such cases can be treated by the analysis in~\cite{maclachlan2007greedy}, they rely on \textit{a posteriori} measurement of the largest $\theta$ for which the constraints in~\cref{eq:linear_program-cgsa} hold, rather than the \textit{a priori} bound that is given by the initial choice of the value of $\theta$ for which we try to solve~\cref{eq:linear_program-cgsa}.  Thus, our goal in
optimizing \cref{eq:linear_program-cgsa} is to improve
complexities of the resulting algorithm (by making the $F$-set larger), subject to the same two-level convergence bound fixed a priori by our choice of $\theta$.  We note that this is further complicated by the complex relationship between two-level and multilevel convergence and complexities.  In particular, when coarsening rates are slow
(corresponding to ``bad'' solutions of~\cref{eq:linear_program-cgsa}),
two-level convergence is generally a bad predictor of multilevel convergence
(since good two-level convergence relies on an assumption of exact solution
of large coarse-grid matrices). We present the supporting numerical
results to explore this connection in \Cref{subsec:struc-grid-algeb-subdom-cgsa}, but note that the work proposed here seeks specifically to improve the quality of solution to~\cref{eq:linear_program-cgsa}, by finding larger $F$-sets that satisfy the constraints for a specified value of $\theta$.

While the convergence bound for AMGr is attractive, it has several
shortcomings.  First and foremost is the strict assumption on diagonal
dominance needed for the convergence guarantee to be valid~---~we
explore cases where this is not the case below in~\cref{subsec: struc
  grid aniso diff-cgsa,subsec:unstructured-cgsa}, and see that these problems
are, indeed, more difficult to handle in this framework.
Second, we emphasize that while the convergence bound obtained
in~\cref{thm:amgr-cgsa} depends only on $\epsilon$ (and, thus,
is independent of mesh size if the coarse grid is generated
by~\cref{alg:greedy_alg-cgsa} or the techniques introduced below with
constant $\theta$). These bounds are generally worse than the observed
convergence for standard AMG approaches.  Using $\theta=0.56$, as we
do for all examples below, gives $\epsilon = 22/3$, and the
convergence bound for $\nu = 1$ in~\cref{eqn:amgr_bound-cgsa} is 0.977.
While this can be improved by using either more relaxation per cycle
or by accelerating convergence with a Krylov method, we consider only
stationary cycles here, accepting poorer convergence factors in
exchange for a direct relationship to the convergence theory
summarized above.  We note that this is also similar to the stationary convergence bound for the method in Napov and Notay~\cite{doi:10.1137/100818509}, which is 0.93 for problems that satisfy the assumptions therein.

%%%%%%%%%%%%%%%%%%%%%%%%%%%%%%%%%%%%%%%%%%%%%%%%%%%%%%%%%%%%%%%%%%%%%%%%%%%%%%

\section{Simulated annealing}\label{sec:sacoarsening-cgsa}

The optimization problem in~\cref{eq:linear_program-cgsa}
is a combinatorial optimization problem.
Because such problems arise in many areas of computational science and
engineering, significant research effort has been devoted to
developing algorithms for their solution, both of general-purpose type
(e.g., branch and bound techniques) and for specific problems (e.g.,
the travelling salesman problem)~\cite{KorteVygen2018}.  For many problems, the size of the
solution space makes exhaustive brute-force algorithms infeasible; for
some such problems, branch and bound techniques may be successful in
paring down the solution space to a more manageable size.
In many cases, however, there are no feasible exact algorithms, and
stochastic search algorithms, such as simulated annealing (SA),
genetic algorithms, or tabu search methods, can be employed~\cite{SchneiderKirkpatrick2006}.
In this work, we apply SA algorithms to approximate the global minimum
of the optimization problem in~\cref{eq:linear_program-cgsa}.

Simulated annealing (SA) is a probabilistic method used to find global
optima of cost functions that might have a large number of local
optima. The SA algorithm randomly generates a state at each iteration
and the cost function is computed for that state. The value of the cost
function for a state determines whether the state is an improvement.
If the current state improves the value of the objective function,
it is accepted to exploit the improved result.
The current state might also be
accepted, with some probability less than one, even if it is worse than
the previous state, however the probability of accepting a bad state decreases
exponentially with the ``badness'' of the state. The purpose of
(sometimes) accepting
inferior states, known as ``exploration'', is to avoid being trapped in local optima, and the inferior
intermediate states are considered in order to give a pathway to a globally better solution.
The total number of iterations of SA depends on an initial ``temperature'' and
the rate of decrease of that temperature. The temperature also affects
the probability of accepting a bad state, with the exploration phase
of the algorithm becoming less probable as the temperature decreases,
to ensure convergence to a global optimum.

%%%%%%%%%%%%%%%%%%%%%%%%%%%%%%%%%%%%%%%%%%%%%%%%%%%%%%%%%%%%%%%%%%%%%%%%%%%%%%

\subsection{Idea and generic algorithm}

Consider the set $F$ of fine points to be the current state.
To find $F$ that maximizes a fitness function $z(F)$, such as $z(F) = |F|$ in~\cref{eq:linear_program-cgsa}, simulated annealing proceeds as shown in~\cref{alg:SA-cgsa}. The temperature $T$ starts at an initial value and decays by a factor $\alpha \in (0,1)$ at each iteration (the ``cooling schedule''). A key choice when using simulated annealing is a method for choosing a neighbor state $\widetilde{F}$ that is close to the current state $F$. In the case of optimizing the choice of a subset $F \subset \Omega$, a straightforward choice is to choose $\widetilde{F}$ by randomly adding or removing an element from $F$. Finally, the function $P(z, \tilde{z}, T)$ is the probability of accepting a move from a current state with fitness $z$ to the new state with fitness $\tilde{z}$. The standard acceptance function is
\begin{equation}\label{eq:SA_prob-cgsa}
  P(z, \tilde{z}, T)\
  = \begin{cases}
    1 & \text{if $\tilde{z} > z$,} \\
    \exp\big(-(z - \tilde{z})/T\big) & \text{otherwise.}
  \end{cases}
\end{equation}
This probability function always accepts transitions that raise the fitness, and sometimes accepts transitions that decrease the fitness.
Occasionally accepting fitness-decreasing transitions is essential to escape local maxima in the energy landscape, with the chance of such transitions being controlled by the current temperature. By analogy with the physical annealing of metals, at high temperatures we will accept almost all fitness-decreasing transitions, allowing exploration of the fitness landscape, but as the temperature cools we will accept fewer of these transitions, until we eventually become trapped near a fitness maximum.
\begin{algorithm}[!ht]
  \caption{Generic simulated annealing (SA) algorithm}\label{alg:SA-cgsa}
  \begin{algorithmic}[1]
    \State Initialize $F$ to a random state and $T$ to an initial temperature
    \For {$n_{\text{steps}}$ iterations}
    \State Randomly pick a neighbor state $\widetilde{F}$ of $F$
    \If {$\text{rand}(0,1) < P(z(F), z(\widetilde{F}), T)$}
    \State $F \gets \widetilde{F}$
      \EndIf
      \State $T \gets \alpha T$
    \EndFor
    \State Return $F$
  \end{algorithmic}
\end{algorithm}

Much work has been devoted to the design of optimal transition probability functions~\cite{Franz2001detSA} and cooling schedules~\cite{deVicente2003thermoSA} in simulated annealing algorithms. While it can be shown that simulated annealing converges to the global maximum with probability one as the cooling time approaches infinity~\cite{Granville1994SAconv}, in practice the performance depends significantly on the selection of the neighbors $\widetilde{F}$ of the current state $F$. A common heuristic for choosing neighbors is to select states with similar fitness, which is more efficient because we are less likely to reject such transitions, although this is in tension with the desire to escape from steep local maxima.

There are also different algorithmic possibilities for the handling of constraints on the state $F$. Given a constraint subset of allowable states, one common method of enforcing the constraint is to pick only neighbor states $\widetilde{F}$ that are in the constraint subset. This method is appropriate if the constraint subset is connected and constraint-satisfying neighbors can always be found, and it is easy to see that this algorithm inherits all theoretical properties of the unconstrained version. An alternative method, which we will use in~\cref{sec:adapt-coars-part-cgsa}, is to allow constraint-violating neighbors to be selected but keep a record of the highest-fitness constraint-satisfying state visited. This method is advantageous when constraint-violating paths in state space allow rapid transitions between (possibly disconnected) areas in the constraint subset. This algorithm also preserves the global guarantees of convergence under the condition that the global maximum is constraint satisfying.

%%%%%%%%%%%%%%%%%%%%%%%%%%%%%%%%%%%%%%%%%%%%%%%%%%%%%%%%%%%%%%%%%%%%%%%%%%%%%%

\subsection{Adaptation to coarse/fine partitioning}
\label{sec:adapt-coars-part-cgsa}

While it is possible to apply SA directly to the optimization problem in~\cref{eq:linear_program-cgsa}, preliminary experiments show that this is inefficient, due to the global coupling of the degrees of freedom in the optimization problem.  To overcome this, we consider a domain decomposition approach to solving the optimization problem, where we first divide the discrete set of degrees of freedom into (non-overlapping) subdomains (the construction of which is considered in detail in the following subsection) and apply SA to the subdomain problems.  These subdomain problems are not independent, therefore information is exchanged about the tentative partitioning on adjacent subdomains; this is accomplished in either an additive (Jacobi-like) or multiplicative (Gauss-Seidel-like) manner.  For sufficiently small subdomains, however, SA can efficiently find near-optimal solutions to localized versions of the optimization in~\cref{eq:linear_program-cgsa}, and we focus on this process below.

SA is run on each subdomain to partition its DoFs into (local) $C$- and $F$-sets.  As the suitability of this partitioning in a global sense naturally depends on decisions being made on adjacent subdomains, we must be careful to consider what happens to DoFs in the global mesh adjacent to those in the subdomain.  If these subdomains already have their own tentative $C/F$ partitions computed (as is expected to be the case on all but the first sweep through the domain), then this partitioning is considered fixed, and used to guide decisions on the current subdomain. If no tentative
partitioning has yet been computed on a neighbouring subdomain, then the points
in the neighbouring subdomain that are connected to some DoFs of the current subdomain are
considered to be $F$-points, while all other DoFs in the neighbouring subdomain are
considered to be $C$-points. These assumptions are used only for the purposes of
computing the partitioning on the current subdomain, to apply hard constraints
on the partitioning.

Specifically, if $\Omega = \{1,2,\ldots,n\}$ is the global set of degrees of
freedom, partitioned into $s$ disjoint subdomains as
$\Omega = \cup_{k=1}^s \Omega_k$, then the optimization problem to be solved on
subdomain $k$ is given as
\begin{align*}
\max_{F_k\subset \Omega_k}& \lvert F_k\rvert, \\
\text{subject to }& |A_{ii}| \geq \theta \sum_{j\in F} |A_{ij}|,\,\forall i \in \bar{F}_k ,
\end{align*}
where information from other subdomains enters in the constraint, both as additional points for which the constraint must be satisfied and in the right-hand side of the constraint where we sum over $j\in F$ (and not $j\in F_k$). Here, $F_k$ is the current set of points in $\Omega_k$ that are in the $F$-set, and we define $C_k = \Omega_k \setminus F_k$ to be the complementary $C$-set. Further, we use $\bar{F}_k$ and $F$ to denote the sets of tentative $F$-points in $\overline{\Omega}_k$ and $\Omega$, respectively, where the standard ``closure'' notation, $\overline{\Omega}_k$, denotes the set of points $j$, such that either $j\in\Omega_k$ or $A_{ij}\neq 0$ for some $i\in\Omega_k$.  In the localized optimization problem above, the set $\bar{F}_k$ contains both the points in $F_k$ and any point $j\in \overline{\Omega}_k\setminus\Omega_k$ such that either $j$ is in the $F_\ell$-set in a neighbouring subdomain, $\Omega_\ell$, that has a tentative partition or $j\in\Omega_\ell$ for a neighbouring subdomain, $\Omega_\ell$, that does not yet have a tentative partition.
We note that a key part of this localization is that we make choices to optimize the partitioning on the local set, $\Omega_k$, but consider the impacts of these choices on the global $F$-set, not only on the local set, $F_k$.  Thus, when we localize to subdomain $\Omega_k$, we consider the constraints on $\bar{F}_k$, including all points in $\Omega$ where the choice of $F_k$ could possibly lead to a constraint violation.

Within an annealing step on a given subdomain, $\Omega_k$, we take the current tentative set $F_k$ as the initial guess for the partitioning, with the exception of the first step, where we take $F_k=\emptyset$ for all $i\in \Omega_k$ (to ensure we start from a configuration that satisfies the constraint). As a benchmark for the annealing process, we initialize $\overline{n}_{\ssf}^{(k)}$ as the largest size of a constraint-satisfying $F$-set seen so far on $\Omega_k$ (taken to be zero on the first iteration), and $z_k$ to be the number of constraint-satisfying $F$-points in the \textit{current} set $\bar{F}_k$. At each annealing step on $\Omega_k$, we swap points in and out of $F_k$, either increasing, decreasing, or maintaining its size, with equal probability.  The basic algorithm for these swaps is given in~\cref{alg:swapFC-cgsa}, where we take the sets $F_k$ and $C_k$ as input, along with values $n_{\ssf}$ and $n_{\ssc}$ giving the numbers of points to swap from $C_k$ to $F_k$ and vice-versa.  We note that there are two possible ways to do this swap, either selecting the elements from $F_k$ and $C_k$ independently and then moving the elements, or first moving the selected elements from $F_k$ and, then, selecting the elements from the updated set $C_k$ to move. We follow the first way as preliminary experiments suggested that it gives slightly better results.
\begin{algorithm}[!ht]
  \caption{swapFC$(F_k,C_k,n_{\ssf},n_{\ssc})$}\label{alg:swapFC-cgsa}
\begin{algorithmic}[1]
  \State $\widetilde{F}_k \gets F_k$
  \State $\widetilde{C}_k \gets C_k$
  \State Randomly select $n_{\ssc}$ points from $F_k$, remove the points from $\widetilde{F}_k$, and add the points to $\widetilde{C}_k$
  \State Randomly select $n_{\ssf}$ points from $C_k$, remove the points from $\widetilde{C}_k$, and add the points to $\widetilde{F}_k$
  \State Return $\widetilde{F}_k$, $\widetilde{C}_k$
\end{algorithmic}
\end{algorithm}

\Cref{alg:anneal-omega-cgsa} shows the complete annealing algorithm, with inputs given by the matrix, $A$, its decomposition into $s$ subdomains, $\{\Omega_k\}_{k=1}^s$, the initial temperature, $T$, and its decay rate, $\alpha$, as well as the number of SA steps to run for each degree of freedom in $A$, $\ndof$ and for each degree of freedom in each cycle, $\ncycle$. In each cycle, annealing is run on every subdomain, with an ordering determined as discussed below. The main annealing step on $\Omega_k$, as given in~\cref{alg:anneal-cgsa}, then takes the form of selecting whether to increase, decrease, or maintain the size of $F_k$, checking if the selected action is possible, and performing it if it is.  To increase the size of $F_k$, we first check that $C_k$ has sufficient entries to move.  If so, we increase the size of $F_k$ by selecting $x+y$ entries of $C_k$ to move to $F_k$ and $y$ entries of $F_k$ to move to $C_k$, for pre-determined values of $x$ and $y$ (typically $x=1$, $y=0$, although these values could also be drawn from a suitable distribution). If the decision is made to swap points, we check that both $F_k$ and $C_k$ have sufficient points to swap, then swap $x$ points from each set into the other (typically $x=1$).  Finally, if the decision is made to decrease the size of $F_k$, we check that it has points to remove, then move $x+y$ points from $F_k$ to $C_k$ and $y$ points from $C_k$ to $F_k$ (ensuring $x+y>0$; typically $x=1$, $y=0$).  Finally, we measure the fitness of the resulting tentative $F$-set and decide whether or not to accept it before decrementing the temperature by the relative factor $\alpha$ and the number of further annealing steps to take by one.
%%%%%%%%%%%%%%%%%%%%
\begin{algorithm}[!ht]
\caption{annealing-on-$\Omega(A, \{\Omega_k\}_{k=1}^s, T, \alpha, \ndof, \ncycle)$}\label{alg:anneal-omega-cgsa}
\begin{algorithmic}[1]
%\State $V_c \gets$ $s$ nodes, randomly selected
%\State $\{\Omega_k\} \gets$ lloyd-aggregation$(A, V_c)$
%\State $\{\Omega_k\} \gets$ partition the DoFs into $s$ subdomains using $A$
\State $F \gets \emptyset$, $F_k \gets \emptyset$ for all $k = 1,\ldots,s$
\State $C_k \gets \Omega_k$ for all $k = 1,\ldots,s$
\State $\overline{n}_{\ssf}^{(k)} \gets 0$, $z_k \gets 0$ for all $k = 1,\ldots,s$
\State $n_{\text{cycle}} \gets {\ndof}/{\ncycle}$ \label{step:cycles}
\For {$i \gets 1, \ldots, n_{\text{cycle}}$}
\For {$k \gets 1, \ldots, s$}
\State $n_{\text{steps}} \gets {\ncycle}\times \lvert \Omega_k \rvert$
\State $F_k, C_k, \overline{n}_{\ssf}^{(k)}, z_k, F, T \gets$ annealing-on-$\Omega_k(F_k, C_k, \overline{n}_{\ssf}^{(k)}, z_k, F, T, \alpha, n_{\text{steps}})$
\EndFor
\EndFor
\State $C \gets \Omega\setminus F$
\State Return $F$, $C$
\end{algorithmic}
\end{algorithm}
%%%%%%%%%%%%%%%%%%%%

\begin{algorithm}[!ht]
\caption{annealing-on-$\Omega_k(F_k, C_k, \overline{n}_{\ssf}^{(k)}, z_k, F, T, \alpha, n_{\text{steps}})$}\label{alg:anneal-cgsa}
\begin{algorithmic}[1]
  %\State Initialize $F_k$ and $C_k$, set $n_{\text{steps}}$, $\alpha$, $T$, $\overline{n}_{\ssf}^{(k)}$, $z$
  %\State Get the global variable $F$
  \For {$n_{\text{steps}}$ iterations}
\State Randomly generate $r \in \{0,1,2\}$ with equal probability
\If {$r=0$ \& $\left|{C}_k\right| \geq x+y$}
\State $\widetilde{F}_k$, $\widetilde{C}_k$ $\gets$ swapFC$(F_k,C_k,x+y,y)$
\EndIf
\If {$r=1$ \& min$(\left|{F}_k\right|, \left|{C}_k\right|) > x$}
\State $\widetilde{F}_k$, $\widetilde{C}_k$ $\gets$ swapFC$(F_k,C_k,x,x)$
\EndIf
\If {$r=2$ \& $\left|{F}_k\right| \geq x+y$}
\State $\widetilde{F}_k$, $\widetilde{C}_k$ $\gets$ swapFC$(F_k,C_k,y,x+y)$
\EndIf
%\State $\bar{F}_k \gets \widetilde{F}_k \cup \big(F\cap (\overline{\Omega}_k\setminus\Omega_k)\big)$
\State Construct tentative $\bar{F}_k$, $\bar{C}_k$ from $\widetilde{F}_k$, $\widetilde{C}_k$, and $F$
%\State $\bar{C}_k \gets \big(\widetilde{C}_k \cup (\overline{\Omega}_k\setminus\Omega_k)\big)\setminus F$
\State $F_k, C_k, \overline{n}_{\ssf}^{(k)}, z, F$ $\gets$ accept$(\Omega_k, \widetilde{F}_k, \widetilde{C}_k, \bar{F}_k, \overline{n}_{\ssf}^{(k)}, z, T, F)$
\State $T \gets \alpha T$
\EndFor
\State Return $F_k$, $C_k$, $\overline{n}_{\ssf}^{(k)}$, $z_k$, $F$, $T$
%\State $F \gets
\end{algorithmic}
\end{algorithm}
%%%%%%%%%%%%%%%%%%%%

The fitness score of a given (tentative) partition over $\overline{\Omega}_k$ is directly calculated as the number of points in the set that satisfy the diagonal dominance criterion, as outlined in~\cref{alg:fitness_score-cgsa}.  In the acceptance step of the algorithm, given as~\cref{alg:accept-cgsa}, we compute the fitness score for the tentative $\bar{F}_k$ and compare it to that of the current (last accepted) set $\bar{F}_k$.  If the fitness score increases or remains same, then we automatically accept the step and update $F_k$, $C_k$, and $z$. In this case, we additionally check if all points in the current $\bar{F}_k$-set are constraint satisfying and if this set increases the value of $\overline{n}_{\ssf}^{(k)}$. If so, we update the value of $\overline{n}_{\ssf}^{(k)}$ (and update the global $F$-set).  If $\tilde{z} < z$, then the step is accepted with a probability that decreases with temperature and $z - \tilde{z}$, but the additional check need not be done.
%%%%%%%%%%%%%
\begin{algorithm}[!ht]
\caption{fitness$(\bar{F}_k)$}\label{alg:fitness_score-cgsa}
\begin{algorithmic}[1]
\State Calculate diagonal dominance for each DoF in the set $\bar{F}_k$
\State $ \tilde{z} \gets $ the number of points that meet the diagonal dominance criterion
\State Return $\tilde{z}$
\end{algorithmic}
\end{algorithm}
%%%%%%%%%%%%%%%
\begin{algorithm}[!ht]
  \caption{accept$(\Omega_k, \widetilde{F}_k, \widetilde{C}_k, \bar{F}_k, \overline{n}_{\ssf}^{(k)}, z, T, F)$}\label{alg:accept-cgsa}
\begin{algorithmic}[1]
\State $\tilde{z}$ $\gets$ fitness$(\bar{F}_k)$
\If {$\tilde{z} \geq z$}
\State $z \gets \tilde{z}$, $F_k \gets \widetilde{F}_k$, $C_k \gets \widetilde{C}_k$
\If {{$\tilde{z}=\left|\bar{F}_k\right|$} $\&$ {$\tilde{z}\geq\overline{n}_{\ssf}^{(k)}$}} \label{line:update_F}
\State $\overline{n}_{\ssf}^{(k)} \gets \tilde{z}$
\State $F \gets (F\setminus \Omega_k)\cup \widetilde{F}_k$
%\State $F_k = \overline{F}_k\cap\Omega_k$
%\State $C_k = \Omega_k\setminus F_k$
%\State $n_{\text{f}}^{(k)} = \left|F_k\right|$
\EndIf
\Else
\State Randomly generate $x\in[0,1]$
\If {$x < {e}^{(-(z - \tilde{z}) / T)}$}
\State $z \gets \tilde{z}$, $F_k \gets \widetilde{F}_k$, $C_k \gets \widetilde{C}_k$
\EndIf
\EndIf
\State Return $F_k, C_k, \overline{n}_{\ssf}^{(k)}, z, F$
\end{algorithmic}
\end{algorithm}
%%%%%%%%%%%

%A key consideration here is how the points in the ``halo'' region of $\overline{\Omega}_k \setminus \Omega_k$ are treated in the algorithm. Within each Gauss-Seidel cycle, we partition the neighbouring subdomains to $\Omega_k$ into those that have already been visited in the current cycle, and those that have not.  In the first Gauss-Seidel cycle, when a neighbouring region of $\Omega_k$ has not yet been visited, we treat all points that are adjacent to points in the halo region as $C$-points, in order to not overly constraint the optimization when we first visit subdomain $\Omega_k$.  In later Gauss-Seidel cycles, these points are partitioned based on the previous Gauss-Seidel cycle, to give a more accurate reflection of the diagonal dominance of the points in the halo region.  For points adjacent to those in the halo in subdomains that have already been visited in the current Gauss-Seidel cycle, we always use the most recent information on the partitioning to determine diagonal dominance of the points in the halo.

%%%%%%%%%%%%%%%%%%%%%%%%%%%%%%%%%%%%%%%%%%%%%%%%%%%%%%%%%%%%%%%%%%%%%%%%%%%%%%

\subsection{Localization and Gauss-Seidel variants}\label{sec:local-gauss-seid-cgsa}

In~\cref{sec:results-cgsa}, we consider both structured-grid and unstructured-grid problems; consequently, we consider both geometric and algebraic decompositions of $\Omega$ into $\{\Omega_k\}_{k=1}^s$. An important advantage of algebraic partitioning, however, is that it can be used to generate deeper multigrid hierarchies, since geometric partitioning can only be used to generate a single coarse level (which has unstructured DoF locations and, thus, does not naturally lead to further geometric partitioning). Here, we outline both strategies.

For structured grids, we can consider geometric decomposition of the fine grid into subdomains.  For problems with (eliminated) Dirichlet boundary conditions, we typically satisfy the diagonal dominance criterion already at all points adjacent to a Dirichlet boundary, so these points are taken as $F$-points from the beginning and not included in the decomposition into subdomains.  A natural strategy is to subdivide the remaining points into square or rectangular subdomains of equal size (modulo boundary/corner cases).  We consider this in~\cref{sec:results-cgsa} for both finite-difference (FD) and bilinear finite-element (FE) discretizations on uniform meshes.  When the number of points (in one dimension) to be decomposed is not evenly divided by the given subdomain size (in one dimension), the right-most and bottom-most subdomains on the mesh are of smaller size, given by the remainder in that division. On structured grids, we can consider either a lexicographical Gauss-Seidel iteration over the subdomains or a four-colored Gauss-Seidel iteration.  In preliminary experiments, the four-colored Gauss-Seidel iteration outperformed the lexicographical strategy by a small margin, so we use this. While the nature of the cycling is not explicitly encoded in the loop over subdomains in~\cref{alg:anneal-omega-cgsa}, we assume that the subdomain indexing in $\{\Omega_k\}$ is consistent with the cycling strategies discussed here.

On both structured and unstructured grids, we also consider algebraic decomposition using Lloyd aggregation to define
the subdomains. Lloyd aggregation, proposed by Bell~\cite{bell2008algebraic} and given in \Cref{alg:lloyd-cgsa}, is a natural application of Lloyd’s algorithm~\cite{lloyd1982least} to subdivide the
DoFs of a matrix into well-shaped subdomains.
Given a desired number of subdomains (or average size per subdomain), the unit-distance graph of the matrix is constructed and one ``center'' point for each subdomain is randomly selected.  Each (tentative) subdomain is then selected as the set of points that are closer to the subdomain center point than to any other subdomain center, using a modified form of the Bellman-Ford algorithm (\Cref{alg:bellFord-cgsa}).  Then, for each subdomain, the center point is reselected (again using a modified form of the Bellman-Ford algorithm), as the current centroid of the subdomain (a DoF having maximum distance from the subdomain boundary, breaking ties arbitrarily).  This process is repeated (reforming subdomains around the new centers, then reassigning centers) until the assignment to subdomains (denoted by the ``membership vector'', $\vec{m}$) has converged.

%%%%%%%%%%%%%%%%%%%%%%%%%%%%%%%%%%%%%%%%%%%%%%

\begin{algorithm}[!ht]
\caption{lloyd-aggregation$(A, V_c)$}\label{alg:lloyd-cgsa}
\begin{algorithmic}[1]
\Repeat
\State $\vec{d}, \vec{m} \leftarrow$ modified-bellman-ford$(A, V_c)$
\State $B \leftarrow \emptyset$
\For {$i,j$ such that $\lvert A_{i,j} \rvert> 0$}
\If {$m_i \ne m_j$}
\State $B \leftarrow B \cup \{i,j\}$
\EndIf
\EndFor
\State $\vec{d},\vec{x} \leftarrow$ modified-bellman-ford$(A, B)$
\State $V_{c} \leftarrow \{i\in \Omega: d_i> d_j \;\;\forall {m_i=m_j}\}$
\Until{no change in $V_{c}$ and $\vec{m}$}
\State Return $\vec{m}$
\end{algorithmic}
\end{algorithm}

%%%%%%%%%%%%%%%%%%%%%%%%%%%%%%%%%%%%%%%%%%%%%%

\begin{algorithm}[!ht]
\caption{modified-bellman-ford$(A, V_c)$}\label{alg:bellFord-cgsa}
\begin{algorithmic}[1]
\State $d_i = \infty$ for all $i = 1,\ldots,\lvert \Omega \rvert$ \Comment{shortest-path distance from node $i$ to nearest center}
\State $m_i = -1$ for all $i = 1,\ldots,\lvert \Omega \rvert$ \Comment{cluster index (membership) containing node $i$}
\For {$c\in V_c$}
\State $d_c \leftarrow 0$
\State $m_c \leftarrow c$
\EndFor
\While {True}
\State done $\leftarrow$ True
\For {$i,j$ such that $\lvert A_{i,j} \rvert> 0$}
\If {$d_i + d_{ij} < d_j$}
\State $d_j \leftarrow d_i + d_{ij}$
\State $m_j \leftarrow m_i$
\State done $\leftarrow$ False
\EndIf
\EndFor
\If {done}
\State Return $\vec{d}, \vec{m}$
\EndIf
\EndWhile
\end{algorithmic}
\end{algorithm}

%%%%%%%%%%%%%%%%%%%%%%%%%%%%%%%%%%%%%%%%%%%%%%

In the algebraic case, a multicoloured Gauss-Seidel iteration strategy is slightly more complicated, since we would need to compute a colouring of the subdomains; hence, we simply use lexicographical Gauss-Seidel to iterate between subdomains, noting that the advantage of the four-color iteration in the structured case is quite small.
While the results in this paper are generated using a serial implementation and lexicographical Gauss-Seidel, the algorithm could easily be parallelized using a multicoloured Gauss-Seidel iteration.

\subsection{Benchmark results}\label{ssec:benchmarks-cgsa}

A natural comparison is with that of the greedy strategy~\cite{maclachlan2007greedy}, which we provide in the following numerical results.  For the structured-grid discretizations, we have also explored optimization ``by hand'', meaning pencil-and-paper analysis of strategies that try and maximize the size of the global F-set while satisfying the constraint.

\begin{figure}[!ht]
\centering
\includegraphics{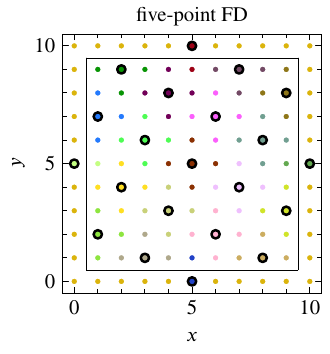}
\qquad
\includegraphics{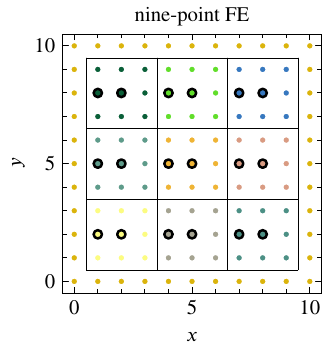}
\caption{Splitting of $C$- and $F$- points for $11\times 11$ meshes from optimization ``by hand''.  At left, X-pentomino coarsening for five-point FD scheme.  At right, $3\times 3$ brick coarsening for nine-point FE scheme. At left, we color by the X-pentominos and at right, we color by the $3\times 3$ bricks. $C$-points are represented by the black circles.}\label{fig:cfnodes_byhand_struc_geom_FDE_11sq}
\end{figure}

For the five-point finite-difference stencil of the Laplacian on a uniform 2D mesh, for any value $\frac{1}{2}<\theta<\frac{4}{7}$, the diagonal dominance
criterion will be satisfied if every $F$-point has at least one $C$-neighbour.
An optimal strategy~\cite{ATag_etal_2021a} for this case arises by dividing the fine mesh into ``X-pentominos'', sets of five grid points with one center point and its four cardinal neighbours, plus edge/corner cases where only a subset of an X-pentomino is needed.  Then, the $C$-set can be selected as the center points of each X-pentomino (or subset of one), with the remaining points assigned to the $F$-set.  In an ideal case (e.g., with periodic boundary conditions, or minimal edge cases), this leads to an $F$-set with size equal to $4/5$ of the size of $\Omega$.  Such a coarsening is shown in the left panel of~\cref{fig:cfnodes_byhand_struc_geom_FDE_11sq}.  Note that some points adjacent to Dirichlet BCs are naturally chosen as coarse points in this strategy, despite their inherent diagonal dominance, because they are center points of an X-pentomino that includes just one point in the region away from the boundary.  These could be equally well treated by making the interior point in these X-pentominos  a $C$-point and leaving the center point on the boundary as an $F$-point, but there is no advantage to doing so.

For the nine-point bilinear finite-element discretization of the Laplacian on a uniform 2D mesh, for any value of $\theta$ less than $4/7$, the diagonal dominance criterion will be satisfied if every $F$-point has at least two $C$-neighbours. Here, we partition the points (again except those adjacent to a Dirichlet boundary) into square subdomains of size $3\times 3$, and select two points consistently in each subdomain as $C$-points.  In an ideal case, where we can fill the domain with $3\times 3$ ``bricks'', this leads to an $F$-set with $7/9$ of the size of $\Omega$.  Such a coarsening is shown in the right panel of~\cref{fig:cfnodes_byhand_struc_geom_FDE_11sq}. When the domain cannot be filled perfectly with $3\times 3$ bricks, the remaining square/rectangular regions still require one or two $C$ points to be selected, reducing the optimal size of the resulting $F$ set.

%%%%%%%%%%%%%%%%%%%%%%%%%%%%%%%%%%%%%%%%%%%%%%%%%%%%%%%%%%%%%%%%%%%%%%%%%%%%%%

\section{Results}\label{sec:results-cgsa}

In the results below, unless stated otherwise, we set the initial (global) temperature at $T=1$, and compute the reduction rate, $\alpha$, so that $\alpha^{n_{\text{ts}}}$ = 0.1, where $n_{\text{ts}}$ is the total number of annealing steps to be attempted.  Thus, these experiments end when the global temperature reaches $0.1$.  In addition, for each problem, we determine $n_\text{ts}$ by fixing a total number of steps per DoF in the system, and report results based on the allotted number of steps per DoF, $\ndof$.  This work is further subdivided into Gauss-Seidel sweeps by fixing values of the number of annealing steps per DoF per sweep, $\ncycle$, with the number of sweeps determined as the ratio of total steps per DoF to steps per DoF per sweep, as on~\cref{step:cycles} of~\cref{alg:anneal-omega-cgsa}.

We first consider structured-grid experiments for both the finite-difference (five-point) and bilinear finite element (nine-point) discretizations of the Laplacian, using a $32\times 32$ mesh to explore how much work is needed to determine near-optimal partitions using both geometric and algebraic divisions into subdomains.  Using these experiments to identify ``best practices'', we then explore how the coarsening algorithm performs as we change the problem size, move from structured to unstructured finite-element meshes, and isotropic to anisotropic operators.

%%%%%%%%%%%%%%%%%%%%%%%%%%%%%%%%%%%%%%%%%%%%%%%%%%%%%%%%%%%%%%%%%%%%%%%%%%%%%%
\subsection{Structured-grid discretizations with geometrically structured subdomains}

\begin{figure}[!ht]
\centering
\includegraphics[width=\textwidth]{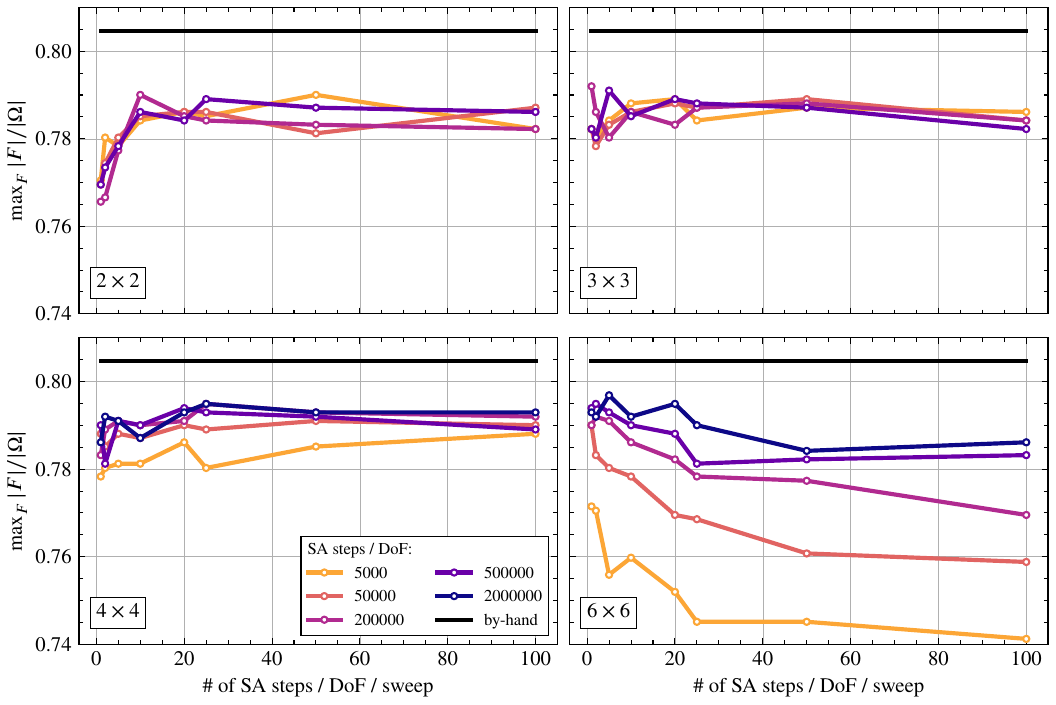}  % chktex 8
\caption{Maximum value of $|F|/|\Omega|$ with number of annealing steps per DoF per GS sweep for different numbers of annealing steps per DoF for the $32\times 32$ uniform-grid five-point finite-difference discretization.  Each panel shows a different size of geometrically chosen subdomain: $2\times 2$ (top-left), $3\times 3$ (top-right), $4\times 4$ (bottom-left), $6\times 6$ (bottom-right).}~\label{fig:maxfptVsastepPsubdomPcyc_FD32_struc_geom_5-50K2-5L2M_2-3-4-6}  % chktex 8
\end{figure}

We start by considering the five-point finite-difference stencil on a fixed (uniform) $32\times 32$ mesh, and consider the effects of changing both the size of the Gauss-Seidel subdomains and the distribution of work in the algorithm.  As a measure of quality of the results, we consider the maximum of the ratio $|F|/|\Omega|$ over all constraint-satisfying $F$-sets generated in a single run of the annealing algorithm.  As a comparison, for this problem, the best optimization ``by hand''  of the size of the $F$-set yields $|F|/|\Omega| = 0.8047$, as depicted by the black lines in \Cref{fig:maxfptVsastepPsubdomPcyc_FD32_struc_geom_5-50K2-5L2M_2-3-4-6}, while the greedy algorithm~\cite{maclachlan2007greedy} yields $|F|/|\Omega| = 0.561$ (not depicted because it is far from the data shown here).  We vary three algorithmic parameters in \Cref{fig:maxfptVsastepPsubdomPcyc_FD32_struc_geom_5-50K2-5L2M_2-3-4-6}, the total number of SA steps per DoF, with values ranging from \num{5000} to \num{2000000}, the number of SA steps per DoF in a single sweep of Gauss-Seidel on each subdomain, and the size of the subdomains used in the Gauss-Seidel sweeps.

  We can draw three conclusions from the data presented in~\Cref{fig:maxfptVsastepPsubdomPcyc_FD32_struc_geom_5-50K2-5L2M_2-3-4-6}.  First, we note that if the subdomains are ``too small'', then there is little benefit in investing substantial work in the SA process, as shown in the top row for $2\times 2$ and $3\times 3$ subdomains.  Here, while there is clearly a small benefit to increasing the number of SA steps per DoF per sweep from 1 to about 10, there is little improvement beyond those results, and little correlation between the quality of partitioning generated and the total amount of work invested.  This occurs consistently in the numerical results throughout this paper: for small subdomain sizes, each subdomain seems to have too little freedom to make adjustments into better global configurations while still satisfying local constraints.  Secondly, when we consider larger subdomains (as in the bottom row), we see that doing more work overall does, indeed, pay off, particularly for the largest subdomains ($6\times 6$, at bottom right).  This is also consistently observed; in particular, that for larger subdomains we see both better overall configurations (if sufficient work is performed) and improvements with larger work budgets.  Finally, for the largest subdomains, we see a clear benefit to doing relatively few SA steps per DoF per cycle.  Thus, in further results, we generally fix the number of SA steps per DoF per cycle to be relatively small, either 1 or 5.

\begin{figure}[!ht]
\centering
\includegraphics[width=\textwidth]{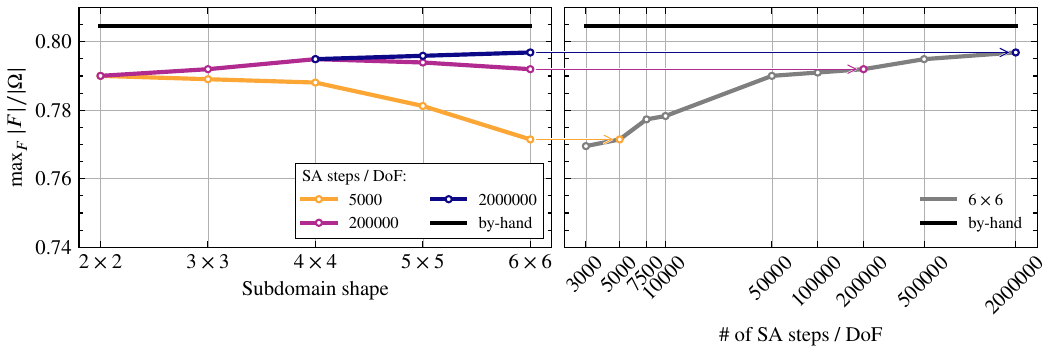}  % chktex 29
\caption{Left: Change in $|F|/|\Omega|$ with subdomain size for the $32\times 32$ uniform-grid five-point finite-difference discretization, using geometric subdomains, with \num{5000}, \num{200000}, and \num{2000000} total SA steps per DoF. Right: Change in maximum number of $F$-points with number of total SA steps per DoF for this problem using $6\times 6$ subdomains.}\label{fig:maxfptVsubdomshape_o2L_maxfptVsastep_5spdpc_6X6_FD32_struc_geom}
\end{figure}

A key question is how to balance the parameters in the SA algorithm to achieve reasonable performance at an acceptable cost.
To examine this, we fix the total number of annealing steps per DoF
and consider the best partitioning achieved by varying the subdomain
size in the left panel
of~\cref{fig:maxfptVsubdomshape_o2L_maxfptVsastep_5spdpc_6X6_FD32_struc_geom}.
Here, we run for a number of different values of the number of SA
steps per DoF per Gauss-Seidel sweep, and take the best partitioning
observed over these grids for each subdomain size (noting that the
optimal choice varies with subdomain size, typically being
larger for small subdomain sizes and smaller for large subdomain
sizes, see~\cref{tab:stepsFig3left} in~\cref{appendix:data} for full details).   When using \num{200000} SA steps per DoF, there is a clear
maximum in the graph for $4\times 4$ subdomains, although the relative
difference in quality in not substantial; however, when using
\num{2000000} SA steps per DoF, we see continued improvement in the
quality of the partitioning up to the $6\times 6$ subdomain case.  In
the right panel
of~\cref{fig:maxfptVsubdomshape_o2L_maxfptVsastep_5spdpc_6X6_FD32_struc_geom},
we look more closely at the convergence of the results with increasing
numbers of SA steps per DoF for the case of $6\times 6$ subdomains,
again taking the best results obtained for different values of the
number of SA steps per DoF per Gauss-Seidel sweep, see~\cref{tab:stepsFig3right} in~\cref{appendix:data} for full details.  Here, we see a
clear improvement in the results up to $\mathcal{O}(10^5)$ SA steps
per DoF, and continued improvement up to \num{2000000} SA steps per
DoF. We recall that we fix the SA temperature reduction rate, $\alpha$,
  so that $\alpha \rightarrow 1$ as the number of SA steps increases,
  yielding the same total reduction in $T$ for each experiment.  Thus,
the results
in~\cref{fig:maxfptVsubdomshape_o2L_maxfptVsastep_5spdpc_6X6_FD32_struc_geom}
are consistent with the expected behaviour of SA, that we can achieve
results arbitrarily close to the global maximizer of our functional
but only if we take many steps and slowly ``cool'' the SA iteration.
For a more practical algorithm, we emphasize the behaviour at
lower numbers of SA steps / DoF, noting that we achieve results
within 5\% of the best-known solution already with only \num{3000} steps /
DoF, and within 2\% at around \num{50000} steps / DoF.
%##
% Left-hand panel of Fig 3 uses:
%subdomains   SAsteps/Dof     SAsteps/DoF/sweep               %subdomains   SAsteps/Dof     SAsteps/DoF/sweep
%-----------------------------------------------              %-----------------------------------------------
% 2\times 2       5000                50                      % 3\times 3       5000                20
%                 50000               100                     %                 50000               50
%                 200000              10                      %                 200000              1
% 4\times 4       5000                100                     % 5\times 5       5000                2
%                 50000               50                      %                 50000               2
%                 200000              25                      %                 200000              5
%                 2000000             25                      %                 2000000             20/50
% 6\times 6       5000                1
%                 50000               1
%                 200000              2
%                 2000000             5
%---------------------------------------------------------------------------------------------------------------
%

\begin{figure}[!ht]
\centering
\includegraphics[width=\textwidth]{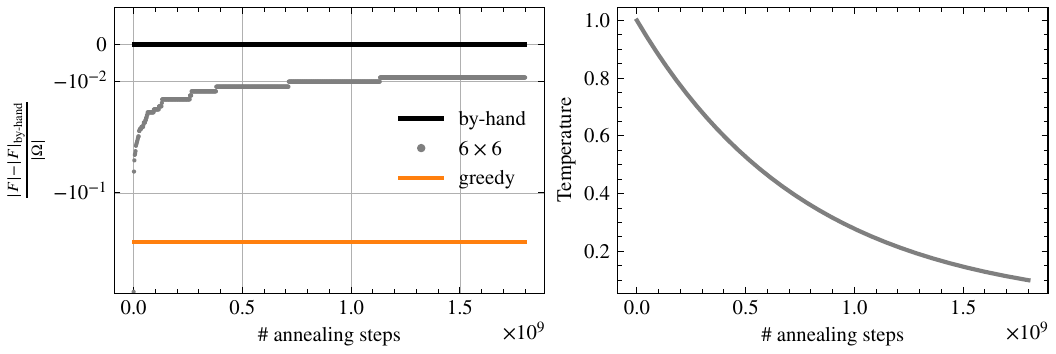} % chktex 29
\caption{Change in $|F|/|\Omega|$ (left) and temperature (right) with number of annealing steps for the $32\times 32$ uniform-grid finite-difference discretization of the Laplacian. At left, the case of SA with five SA steps per DoF per GS sweep using $6 \times 6$ subdomain size is shown, along with greedy as a baseline. Note that the vertical axis on the left plot uses a mixed log-linear scale for clarity, with a break at $-10^{-2}$.}\label{fig:fpts-tempVSsastep_FD32_geom_6X6_2M_5_symlog}
\end{figure}

A natural question that arises from the right-hand panel of~\Cref{fig:maxfptVsubdomshape_o2L_maxfptVsastep_5spdpc_6X6_FD32_struc_geom} is whether the best meshes obtained occur early or late in the annealing process.  That is, while we clearly see benefit from slow ``cooling'' of the annealing process (many SA steps per DoF), it is important to identify \textit{when} the optimal results are obtained during the annealing process.  \Cref{fig:fpts-tempVSsastep_FD32_geom_6X6_2M_5_symlog} shows the annealing history for a sample run, on the $32\times 32$ uniform grid five-point finite-difference stencil, with \num{2000000} SA steps per DoF (for a total of almost $2\times 10^9$ annealing steps).  At right, we see the temperature decay, following an exponential curve from $T=1.0$ at the first step to $T=0.1$ at the final step.  At left, we plot the changes in the ratio $|F|/|\Omega|$ compared to the optimized by hand mesh for this grid (also shown).  For comparison, the ratio from the greedy algorithm is also shown.  We see that after an initial rapid improvement in the value of the ratio, there is a secondary period where performance improves notably but steadily, up to between $5\times 10^8$ and $10^9$ total annealing steps.  This demonstrates that many annealing steps are, indeed, needed to reach the best partitioning seen here, although a good partitioning is still found after many fewer steps.

\begin{figure}[!ht]
\centering
\includegraphics[width=\textwidth]{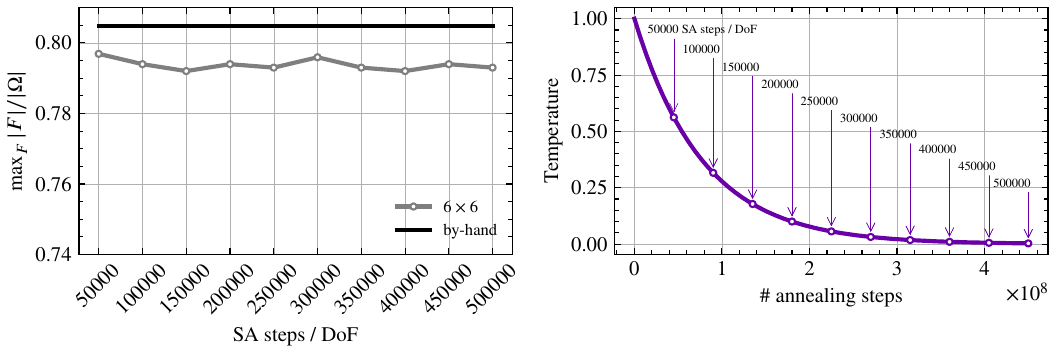}
\caption{Change in $|F|/|\Omega|$ with SA steps per DoF for a fixed temperature decay rate (left) and temperature (right) for the $32\times 32$ uniform-grid finite-difference discretization of the Laplacian.  Here, one SA step per DoF per sweep is used, and the temperature decay rate is fixed, with $\alpha = 0.1^{(1.0/(200000\times32\times32))}$. The circles in the right figure show the stopping temperatures for the annotated SA steps per DoF.}\label{fig:fpts_vs_sasteps_per_dof_6X6_FD32_struc_geom_fixed_temp_decay_rate}
\end{figure}

Thus far, we have only considered the choice of $\alpha$ prescribed at the beginning of this section, with the decay rate chosen to yield a fixed temperature decay by a factor of 10 over the total number of SA steps given.  In contrast, \cref{fig:fpts_vs_sasteps_per_dof_6X6_FD32_struc_geom_fixed_temp_decay_rate} considers fixing the decay rate to be that used for \num{200000} SA steps per DoF, $\alpha = 0.1^{(1.0/(200000\times32\times32))}$, but then running between 4 times fewer and 2.5 times more total SA steps, with one SA step per DoF per cycle, yielding final temperature values between $0.5623$ and $0.0032$.  At left of~\cref{fig:fpts_vs_sasteps_per_dof_6X6_FD32_struc_geom_fixed_temp_decay_rate}, we observe little correlation between the total number of SA steps per DoF and the resulting value of $|F|/|\Omega|$, with all values between $0.79$ and $0.80$, comparable to those seen for similar SA budgets in~\cref{fig:maxfptVsubdomshape_o2L_maxfptVsastep_5spdpc_6X6_FD32_struc_geom}.  While this appears to offer an opportunity for saving some cost in the SA algorithm (by using fewer Gauss-Seidel sweeps to get comparable results), we note that using \num{50000} SA steps per DoF is already prohibitively expensive for an ``online cost'' for this algorithm, requiring more than 30 minutes to compute the partitioning for this relatively small problem.
%%%%%%%%%%%%%%%%%%%%%%%%%%%%%%%%%%%%%%%%%%%%%%%%%%%%%%%%%%%%%%%%%%%%%%%%%%%%%%
%

%
%\begin{figure}[!ht]
%\centering
%\includegraphics[width=\textwidth]{fpts_vs_temp_decay_rate_6X6_FD32_struc_geom_fixed_sastep.pdf}
%\caption{Change in $|F|/|\Omega|$ with temperature decay rate (left) and temperature profiles for the $32\times 32$ uniform-grid finite-difference discretization of the Laplacian. At left, the case of SA with one SA step per DoF per GS sweep and \num{200000} SA steps per DoF is used.}\label{fig:fpts_vs_temp_decay_rate_6X6_FD32_struc_geom_fixed_sastep}
%\end{figure}
%%%%%%%%%%%%%%%%%%%%%%%%%%%%%%%%%%%%%%%%%%%%%%%%%%%%%%%%%%%%%%%%%%%%%%%%%%%%%%

\begin{figure}[!ht]
\centering
\includegraphics[width=0.45\textwidth]{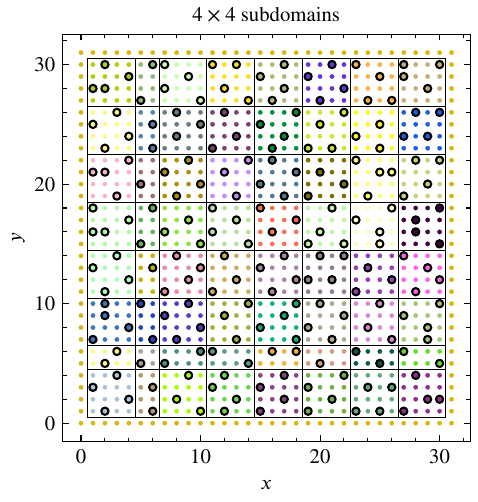} % chktex 8 chktex 29
\qquad
\includegraphics[width=0.45\textwidth]{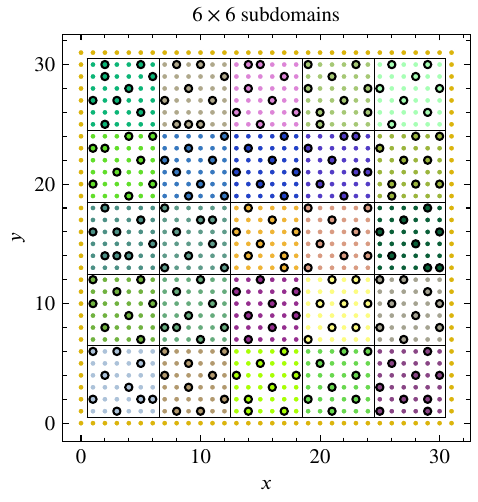} % chktex 8 chktex 29
\caption{Grid partitioning for $32\times 32$ uniform grid with five-point FD stencil using two million SA steps per DoF. At left, the partitioning is generated using primarily $4\times 4$ subdomains, with 25 SA steps per DoF per cycle.  At right, the partitioning is generated using $6\times 6$ subdomains, with 5 SA steps per DoF per cycle.  The grid at left has 814 $F$-points, while that at right has 816.}\label{fig:cfnodes_struc_32X32_FD_geom_4-6_2000000_25-5}
\end{figure}

We next address the nature of the grids generated by annealing, and whether
they resemble grids that could be selected geometrically for this problem.  \Cref{fig:cfnodes_struc_32X32_FD_geom_4-6_2000000_25-5} shows two of the grids generated, along with the subdomains used in their generation.  These represent the ``best'' grids found by the annealing procedure, with ratios of $|F|/|\Omega|$ of around $0.795$.  While these grids yield competitive ratios, there is no clear global geometric pattern, nor obvious relationship to the best  optimized by hand grid shown in~\cref{fig:cfnodes_byhand_struc_geom_FDE_11sq}.  Furthermore, there are no clear improvements of these grids that could be readily made, such as single coarse points that could be omitted without leading to constraint violations.  This suggests that the energy landscape for this problem is likely dominated by locally optimal configurations that are separated by states with constraint violations and/or sharp changes in energy.

\begin{figure}[!ht]
\centering
\includegraphics[width=\textwidth]{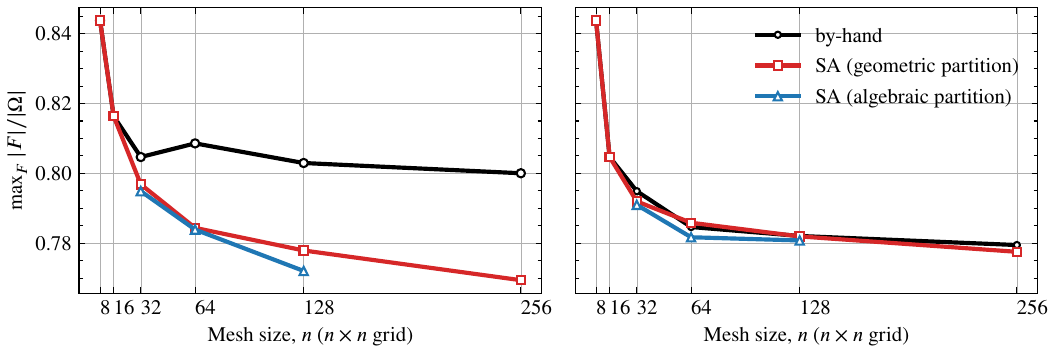}
\caption{Change in $|F|/|\Omega|$ with mesh size for FD (left) and FE (right) discretizations.}\label{fig:maxfptVmeshsize_FDE_struc_geom_lloyd3}
\end{figure}

\Cref{fig:maxfptVmeshsize_FDE_struc_geom_lloyd3} shows how the performance of the annealing algorithm scales with problem size, for both the geometric choice of subdomains for the finite-difference discretization discussed so far, and for the finite-element discretization.  In addition, an algebraic selection (discussed below) is added for comparison.  All methods (including the optimization by hand) perform relatively well for small meshes, but degrade as the mesh size increases.  The amount of work needed to achieve these results with the annealing algorithm also increases with grid size.  For the $8\times 8$ mesh, using a single (global) subdomain, equally-good partitions to the optimization by hand can be found with just \num{2000} annealing steps per DoF (and fewer when more subdomains are used).  For the $16\times 16$ mesh, more work and larger subdomains are needed to achieve such performance.  As noted above, with small subdomains even using \num{200000} annealing steps per DoF does not achieve performance equal to the by-hand partitioning on the $16\times 16$ mesh.  For larger subdomains, the algorithm does equal the results of the by-hand partitioning, but with increasing work as subdomain size increases: for $4\times 4$ subdomains, \num{10000} annealing steps per DoF are needed, while \num{50000} annealing steps per DoF are needed for $5\times 5$ subdomains.  For $6\times 6$ subdomains, even using \num{200000} annealing steps per DoF, we could not recover results matching the by-hand partitioning, although we speculate that this would have occurred with even more work invested.  For larger domain sizes, the results shown in~\Cref{fig:maxfptVmeshsize_FDE_struc_geom_lloyd3} represent the best results found for a given grid over all runs with varying subdomain sizes, total SA steps per DoF, and SA steps per DoF per GS sweep.
While these best results are, in general, achieved with the
largest allocations of SA steps per DoF tried in our experiments,
the marginal benefit of considering such large amounts of work to
generate the coarsenings are quite low.  Here, in all cases where we
invested the computational effort to explore the question, we found
less than a 1\% improvement in $|F|/|\Omega|$ when increasing beyond
$\mathcal{O}(10^5)$ SA steps per DoF (and, in many cases, the
improvement was only by one or two $F$-points).

\begin{figure}[!ht]
\centering
\includegraphics[width=\textwidth]{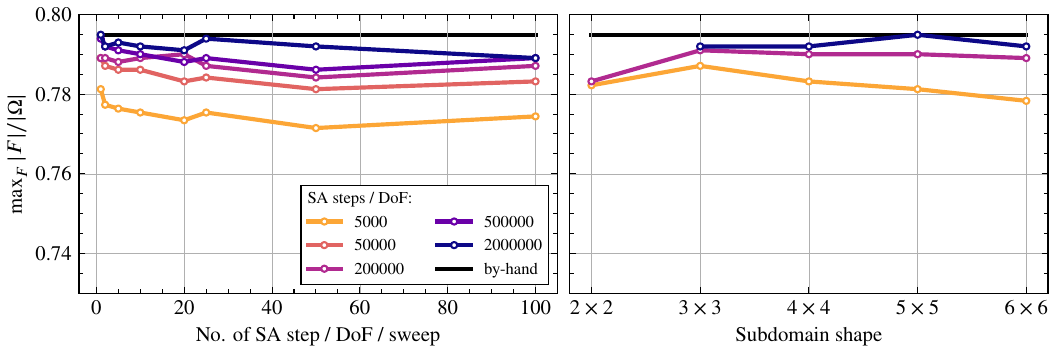}
\caption{Quality of coarsening for the $32\times 32$ uniform-grid nine-point finite-element discretization, using geometric subdomains. Left: Maximum value of $|F|/|\Omega|$ with number of annealing steps per DoF per GS sweep for different numbers of annealing steps per DoF using $5\times 5$ geometric subdomains. Right: Change in $|F|/|\Omega|$ with subdomain size and total number of SA steps per DoF.}\label{fig:maxfptVsastepPsubdomPcyc_FE32_5_maxfptVsubdomshape_struc_geom}
\end{figure}

\begin{figure}[!ht]
\centering
\includegraphics[width=0.45\textwidth]{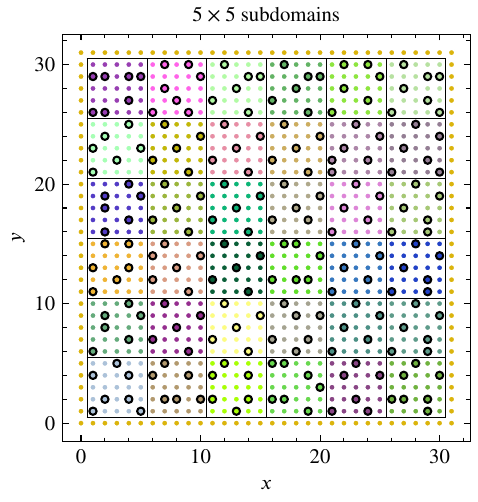}  % chktex 29 chktex 8
\qquad
\includegraphics[width=0.45\textwidth]{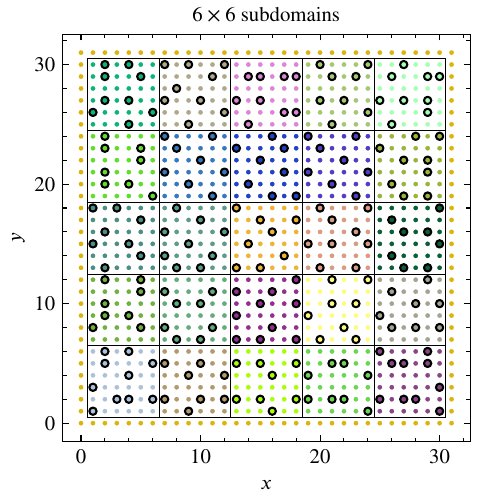}  % chktex 29 chktex 8
\caption{Grid partitioning for $32\times 32$ uniform grid with nine-point FE stencil using two million SA steps per DoF. At left, the partitioning is generated using $5\times 5$ subdomains, with one SA step per DoF per cycle.  At right, the partitioning is generated using $6\times 6$ subdomains, also with one SA step per DoF per cycle.  The grid at left has 814 $F$-points, while that at right has 811.}\label{fig:cfnodes_struc_32X32_FE_geom_5-6_2000000_1}
\end{figure}

For the nine-point finite-element stencil, we see similar results to those reported above and, consequently, do not include figures detailing the individual experiments in as much detail.  Most notably, the best partitioning that we achieve by hand is a slightly worse in this case (as detailed in~\Cref{ssec:benchmarks-cgsa}), and the partitioning using the greedy algorithm is better than in the FD case, yielding $|F|/|\Omega| = 0.752$.  At left of \Cref{fig:maxfptVsastepPsubdomPcyc_FE32_5_maxfptVsubdomshape_struc_geom}, we show an analogous figure to \Cref{fig:maxfptVsastepPsubdomPcyc_FD32_struc_geom_5-50K2-5L2M_2-3-4-6} for the case of $5\times 5$ subdomains.  Here, we observe the same stratification. However, we are able to recover the same quality of partitioning as the best by-hand partitioning using $5\times 5$ subdomains and \num{2000000} SA steps per DoF.  Also note that we see the same mild dependence on the number of SA steps per DoF per sweep as we did in the FD case.  The right panel of \Cref{fig:maxfptVsastepPsubdomPcyc_FE32_5_maxfptVsubdomshape_struc_geom} shows how the quality of partitioning changes with subdomain size and total work budget, with SA steps per DoF per sweep reported in~\cref{tab:stepsFig8} in~\cref{appendix:data}.  Similarly to the FD case, there is an improvement in performance with additional work for larger subdomain sizes, but that improvement stagnates with additional work for smaller subdomain sizes.  Sample grids generated for the finite-element case, including colouring to indicate subdomain choice, are shown in~\Cref{fig:cfnodes_struc_32X32_FE_geom_5-6_2000000_1}.
%##
% the following SA steps / DoF / sweep are used:
%subdomains   SAsteps/Dof     SAsteps/DoF/sweep               %subdomains   SAsteps/Dof     SAsteps/DoF/sweep
%-----------------------------------------------              %-----------------------------------------------
% 2\times 2       5000                20                      % 3\times 3       5000                10
%                 50000               25                      %                 50000               25
%                 200000              10/20                   %                 200000              2
                                                              %.                2000000             50
% 4\times 4       5000                5/20                    % 5\times 5       5000                1
%                 50000               5/20                    %                 50000               1
%                 200000              50                      %                 200000              20
%                 2000000             50                      %                 2000000             1
% 6\times 6       5000                1
%                 50000               1
%                 200000              2
%                 2000000             1
%---------------------------------------------------------------------------------------------------------------
%

Finally, we verify that the result two-level AMGr algorithms are at least as effective as predicted in theory, by measuring asymptotic convergence factors ($\rho$), grid complexities ($C_{\text{grid}}$, equal to the ratio of sum of the number of DoFs on each level of the hierarchy to that on the finest level), and operator complexities
($C_{\text{op}}$, equal to the ratio of the sum of the number of nonzero entries in the system matrix on each level of the hierarchy to that on the finest level).  We approximate $\rho$ by running the cycle with a random initial guess, $x^{(0)}$, and zero right-hand side, and then compute
$\left(\left\| x^{(k)}\right\|/\left\| x^{(0)}\right\|\right)^{1/k}$, where $x^{(k)}$ is the approximation (to the true solution, which is the zero vector) after $k$ cycles.  Because the convergence factors of AMGr are somewhat larger than those for classical AMG, we take $k=800$ to ensure that we sample the asymptotic behaviour suitably.
\Cref{tab:conv_isot_struc_geom-lloyd_2lev-cgsa} reports this data for two-level cycles for both the FD and FE discretizations, for both the geometric subdomain choice considered here and the algebraic subdomain choice discussed next.  Considering the geometric subdomain choice, we see grid-independent convergence factors of about 0.9 for the FD case and 0.7 for the FE case.  While these are notably worse than are observed for typical multigrid methods for these problems, they are consistent with the existing results for AMGr and, in particular, conform with the convergence rate bound from~\cref{thm:amgr-cgsa} of 0.977 for $\theta=0.56$ and $\nu = 1$. Notably, neither the convergence factor nor the measured complexities degrade substantially with problem size. We note that, in subsequent work, we have improved convergence of AMGr for these and other model problems~\cite{TZaman_etal_2022a}.
\begin{table}[!ht]
	\caption{Performance of two-level AMGr on test matrices from discretizations of the 2D Laplacian.}
	\centering
	%\begin{adjustbox}{width=0.6\textwidth}
	    %\resizebox{\textwidth}{!}{\begin{tabular}{|c|*{8}{c|}}
	    %\small\begin{tabular*}{\textwidth}{l@{\extracolsep{\fill}}*{9}{c}}
		{\begin{tabular}{c c |c c c |c c c }
			\toprule
			 &  & \multicolumn{3}{c|}{Geometric Subdomains} & \multicolumn{3}{c}{Algebraic Subdomains} \\
      Scheme & Grid size   & $\rho$ & $C_{\text{grid}}$ & $C_{\text{op}}$ & $\rho$ & $C_{\text{grid}}$ & $C_{\text{op}}$ \\
			\midrule
			   & $8\times8$ & 0.85 & 1.16 & 1.12 & 0.85 & 1.16 & 1.12 \\
			   & $16\times16$ & 0.86 & 1.19 & 1.18 & 0.86 & 1.19 & 1.18 \\
			  FD & $32\times32$ & 0.88 & 1.20 & 1.20 & 0.88 & 1.21 & 1.20 \\
			   & $64\times64$ & 0.89 & 1.22 & 1.22 & 0.88 & 1.22 & 1.22 \\
			   & $128\times128$ & 0.89 & 1.22 & 1.23 & 0.88 & 1.22 & 1.23 \\
			   %& $256\times256$ & 0.881 & 1.239 & 1.251 &  &  & \\
         \midrule
			   & $8\times8$ & 0.70 & 1.16 & 1.15 & 0.70 & 1.16 & 1.15 \\
			   & $16\times16$ & 0.63 & 1.20 & 1.22 & 0.63 & 1.20 & 1.22 \\
			  FE & $32\times32$ & 0.67 & 1.21 & 1.23 & 0.69 & 1.21 & 1.25 \\
			   & $64\times64$ & 0.71 & 1.21 & 1.25 & 0.72 & 1.22 & 1.27 \\
			   & $128\times128$ & 0.71 & 1.22 & 1.27 & 0.71 & 1.22 & 1.27 \\
			   %& $256\times256$ & 0.715 & 1.224 & 1.274 &  &  & \\
			\bottomrule
		\end{tabular}}\label{tab:conv_isot_struc_geom-lloyd_2lev-cgsa}
\end{table}
%%%%%%%%%%%%%%%%%%%%%%%%%%%%%%%%%%%%%%%%%%%%%%%%%%%%%%%%%%%%%%%%%%%%%%%%%%%%%%

\subsection{Structured-grid discretizations with algebraically chosen subdomains}\label{subsec:struc-grid-algeb-subdom-cgsa}

\begin{figure}[!ht]
\centering
\includegraphics[width=\textwidth]{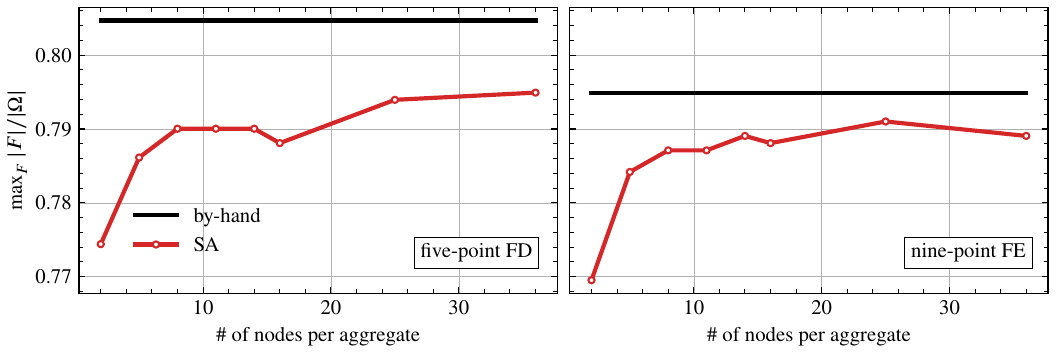}
\caption{Change in maximum number of $F$-points with subdomain size for algebraic subdomain selection on $32\times 32$ meshes. Left and right
  figures are for FD and FE schemes, respectively, showing the largest value of $|F|/|\Omega|$ attained over experiments with fixed subdomain size and varying the total number of SA steps per DoF and SA steps per sweep, as in the geometric subdomain case.}\label{fig:maxfptVsubdomshape-nnpagg_FDE32_struc_lloyd}
\end{figure}

We next consider partitioning the structured-grid problems using an algebraic choice of the subdomains based on Lloyd aggregation. \Cref{fig:maxfptVsubdomshape-nnpagg_FDE32_struc_lloyd} shows the change in the maximum number of $F$-points (scaled by the total number of
DoFs) with the change of the subdomain size for both the FD and FE discretizations (left and right, respectively).  Similarly to the case of geometric partitioning, we see relatively poor performance for small subdomain sizes, regardless of the work allocated to the SA process.  For larger subdomain sizes, we see improving results with number of SA steps per DoF, as seen above.  Also as seen above, the optimal subdomain size varies with total work allocation, increasing as we increase the amount of work per DoF, but even with \num{2000000} SA steps per DoF, we do not see the best performance with largest subdomains for the FE discretization.  Considering variation in problem size, \Cref{fig:maxfptVmeshsize_FDE_struc_geom_lloyd3} shows that, as in the geometric subdomain case, we see some decrease in performance as problem size grows, but that this decrease seems to (mostly) plateau at larger problem sizes.  In comparison with the geometric subdomain choice, we see some small degradation in performance with algebraically chosen subdomains, particularly in the FD case, but it is small in comparison with the optimality gap between the optimization by hand solutions and those generated by SA\@.

Performance of the resulting two-grid cycles is tabulated in \Cref{tab:conv_isot_struc_geom-lloyd_2lev-cgsa}, where we see very comparable convergence factors, as well as grid and operator complexities, as in the geometric subdomain case.
\Cref{tab:conv_isot_struc_lloyd_3lev_V_W-cgsa} tabulates convergence factors for three-level V- and W-cycles (denoted $\rho_V$ and $\rho_W$, respectively), along with three-level grid and operator
  complexities. For the FD case, we see some degradation in convergence when using V-cycles, while W-cycles give convergence factors similar to those in the two-level case.  For the FE problem, there is less degradation for V-cycles, but still W-cycles are required to get convergence factors comparable to the two-level case.  As we have increased the depth of the multigrid hierarchy, the grid and operator complexities in \Cref{tab:conv_isot_struc_lloyd_3lev_V_W-cgsa} are naturally larger than those in \Cref{tab:conv_isot_struc_geom-lloyd_2lev-cgsa}, but the growth in these complexity measures is consistent with optimally scaling multigrid methods. Performance of multilevel V- and W-cycles is also shown in \Cref{tab:conv_isot_struc_lloyd_3lev_V_W-cgsa}, where coarsening is continued until the number of nodes falls below \num{100}, yielding a four-level cycle for the $32\times 32$ grid, five-level cycle for $64\times 64$, and six-level cycle for $128\times 128$.  The resulting convergence factors for both FD and FE cases remain similar to those for three-level cycles, while complexity measures grow consistently with the number of levels.
\begin{table}[!ht]
	\caption{Performance of three-level and multilevel AMGr on test matrices from discretizations of the 2D Laplacian.}
	\centering
	%\begin{adjustbox}{width=0.6\textwidth}
	    %\resizebox{\textwidth}{!}{\begin{tabular}{|c|*{8}{c|}}
	    %\small\begin{tabular*}{\textwidth}{l@{\extracolsep{\fill}}*{9}{c}}
		{\begin{tabular}{c c|c c c c|c c c c}
			\toprule
			       &           & \multicolumn{4}{c|}{Three-level cycles} & \multicolumn{4}{c}{Multilevel cycles} \\
      Scheme & Grid size & $\rho_V$ & $\rho_W$ & $C_{\text{grid}}$ & $C_{\text{op}}$ & $\rho_V$ & $\rho_W$ & $C_{\text{grid}}$ & $C_{\text{op}}$ \\
			\midrule
			   & $32\times32$ & 0.92 & 0.90 & 1.25 & 1.26 & 0.93 & 0.90 & 1.26 & 1.27 \\
			FD & $64\times64$ & 0.93 & 0.91 & 1.27 & 1.28  & 0.94 & 0.91 & 1.28 & 1.31 \\
			   & $128\times128$ & 0.94 & 0.91 & 1.28 & 1.30  & 0.95 & 0.92 & 1.30 & 1.35 \\
			   %& $256\times256$ & 0.916 & 1.301 & 1.338 & 0.880 & 1.301 & 1.338 \\
			\midrule
			   & $32\times32$ & 0.76 & 0.72 & 1.25 & 1.31  & 0.76 & 0.72 & 1.26 & 1.31 \\
			FE & $64\times64$ & 0.76 & 0.73 & 1.27 & 1.35  & 0.77 & 0.73 & 1.28 & 1.37 \\
			   & $128\times128$ & 0.79 & 0.75 & 1.28 & 1.37 & 0.79 & 0.75 & 1.29 & 1.41 \\
			   %& $256\times256$ &  &  &  &  &  & \\
			\bottomrule
		\end{tabular}}\label{tab:conv_isot_struc_lloyd_3lev_V_W-cgsa}
\end{table}

For the multilevel results in~\cref{tab:conv_isot_struc_lloyd_3lev_V_W-cgsa}, we use \num{200000} SA steps per DoF in all cases, except on the finest meshes for the $32\times 32$ and $64\times 64$ grids, where \num{2000000} SA steps per DoF were used.  (Similar results are also seen without these ``extra'' steps on these problems, however.)  For the $64\times 64$ mesh, using \num{200000} SA steps per DoF on all levels takes about 18 hours of ``offline'' time to compute the coarse meshes, in comparison to ``online'' solution times of only 0.003 seconds for a V-cycle and 0.011 seconds for a W-cycle. \Cref{tab:conv_isot_struc_lloyd_multilev_V_W_200saperdof-cgsa} presents results for the same experiment, but using only \num{200} and \num{2000} SA steps per DoF, to investigate how performance is changed when using ``bad'' solutions to the optimization problem in~\cref{eq:linear_program-cgsa}.  Here, we see slight improvements in convergence factors in comparison to those in~\cref{tab:conv_isot_struc_lloyd_3lev_V_W-cgsa}; however, using fewer annealing steps leads to much larger grid and operator complexities than observed above.  We note here that using fewer annealing steps also leads to results that are more heavily influenced by the random nature of the SA process; here, we report complexities from runs used to generate W-cycle convergence factors, which vary slightly from those of an independent run to generate V-cycle convergence factors, but by no more than $0.01$ in $C_{\text{grid}}$ and $0.05$ in $C_{\text{op}}$.

%%%%%%%%%%%%%%%%%%%%%%%%%%%%%%%%%%%%%%%%%%%%%%%%%%%%%%%%%%%%%%%%%%%%%%%%%%%%%%
%
\begin{table}[!ht]
	\caption{Performance of multilevel AMGr on test matrices from discretizations of the 2D Laplacian using \num{200} and \num{2000} SA steps per DoF.}
	\centering
		{\begin{tabular}{c c|c c c c|c c c c}
			\toprule
			       &           & \multicolumn{4}{c|}{\num{200} SA steps per DoF} & \multicolumn{4}{c}{\num{2000} SA steps per DoF} \\
      Scheme & Grid size & $\rho_V$ & $\rho_W$ & $C_{\text{grid}}$ & $C_{\text{op}}$ & $\rho_V$ & $\rho_W$ & $C_{\text{grid}}$ & $C_{\text{op}}$\\
			\midrule
			   & $32\times32$ & 0.89 & 0.86 & 1.42 & 1.62 & 0.91 & 0.88 & 1.33 & 1.44 \\
			FD & $64\times64$ & 0.89 & 0.88 & 1.46 & 1.86 & 0.92 & 0.88 & 1.37 & 1.51 \\
			   & $128\times128$ & 0.91 & 0.88 & 1.48 & 2.08 & 0.93 & 0.89 & 1.38 & 1.58 \\
			\midrule
			   & $32\times32$ & 0.72 & 0.68 & 1.34 & 1.52 & 0.73 & 0.71 & 1.29 & 1.40 \\
			FE & $64\times64$ & 0.74 & 0.70 & 1.37 & 1.65 & 0.76 & 0.70 & 1.32 & 1.49 \\
			   & $128\times128$ & 0.77 & 0.70 & 1.38 & 1.76 & 0.76 & 0.72 & 1.33 & 1.56 \\
			\bottomrule
		\end{tabular}}\label{tab:conv_isot_struc_lloyd_multilev_V_W_200saperdof-cgsa}
\end{table}

\begin{figure}[!ht]
\centering
\includegraphics{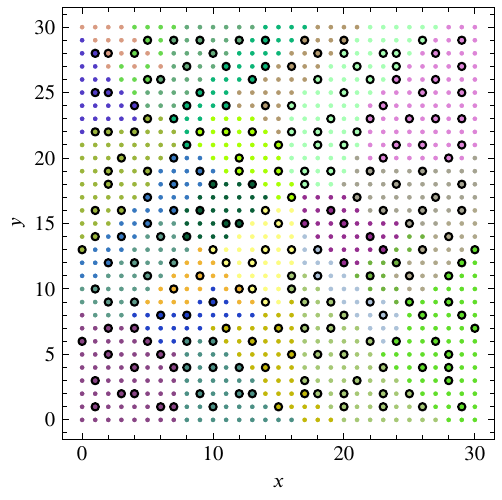}
\caption{Partitioning for the isotropic problem on a $31\times 31$ uniform grid with jumping coefficients, generated using \num{200000} SA steps per DoF, 1 SA step per DoF per GS sweep, and subdomains with an average of 36 points per subdomain. $C$-points are represented by the black circles.}\label{fig:cfnodes_jump_coeff_struc_iso_33X33_FE_lloyd_36_26_200000_1}
\end{figure}

For our final isotropic, structured-grid test problem, we consider the FE discretization of the two-dimensional isotropic diffusion problem, $-\nabla \cdot \boldsymbol{K}(x, y)\nabla u(x, y) = f(x,y)$, in the
  domain $[0,1]\times [0,1]$ with Dirichlet boundary conditions and piecewise constant (``jumping'') coefficient, $\boldsymbol{K}(x,y)$.  To determine $\boldsymbol{K}(x,y)$, we select 20\% of the elements at random, and set $\boldsymbol{K}(x,y) = 10^{-8}$ in these elements, with value $1$ in the remaining 80\% of the elements, matching one of the test problems from~\cite{maclachlan2007greedy}. We partition the nodes using \num{200000} SA steps per DoF, and show the first coarse mesh chosen using $1$ SA step per DoF per sweep for the $31\times 31$ grid in~\cref{fig:cfnodes_jump_coeff_struc_iso_33X33_FE_lloyd_36_26_200000_1}. Convergence factors and grid and operator
complexities for multilevel V- and W-cycles are shown in~\cref{tab:conv_isot_struc_lloyd_jump_coeff_multilev_V_W-cgsa}.  We note, in particular, that these results are quite comparable to those shown in~\cref{tab:conv_isot_struc_lloyd_3lev_V_W-cgsa} for the case of $\boldsymbol{K}(x,y) = 1$ everywhere.  In comparison to the results presented in~\cite{maclachlan2007greedy}, we see larger convergence factors here (0.78 for the $127\times 127$ grid, in comparison to 0.59 reported there), but with lower complexities ($C_{\text{op}} = 1.42$ here, compared to $1.75$ there).

%%%%%%%%%%%%%%%%%%%%%%%%%%%%%%%%%%%%%%%%%%%%%%%%%%%%%%%%%%%%%%%%%%%%%%%%%%%%%%
%%%%%%%%%%%%%%%%%%%%%%%%%%%%%%%%%%%%%%%%%%%%%%%%%%%%%%%%%%%%%%%%%%%%%%%%%%%%%%
%
\begin{table}[!ht]
	\caption{Performance of multilevel AMGr on test matrices from discretizations of the 2D Laplacian with jumping coefficients.}
	\centering
		{\begin{tabular}{c c|c c c c}
			\toprule
      Scheme & Grid size & $\rho_V$ & $\rho_W$ & $C_{\text{grid}}$ & $C_{\text{op}}$ \\
			\midrule
			   & $31\times31$ & 0.73 & 0.71 & 1.27 & 1.33 \\
			FE & $63\times63$ & 0.76 & 0.73 & 1.29 & 1.38 \\
			   & $127\times127$ & 0.78 & 0.76 & 1.29 & 1.42 \\
			   %& $256\times256$ &  &  &  &  &  &  &  &  \\
			\bottomrule
		\end{tabular}}\label{tab:conv_isot_struc_lloyd_jump_coeff_multilev_V_W-cgsa}
\end{table}

%%%%%%%%%%%%%%%%%%%%%%%%%%%%%%%%%%%%%%%%%%%%%%%%%%%%%%%%%%%%%%%%%%%%%%%%%%%%%%

\subsection{Anisotropic problems on structured grids}\label{subsec: struc grid aniso diff-cgsa}

Next, we consider the two-dimensional anisotropic diffusion problem, $-\nabla \cdot \boldsymbol{K}(x, y)\nabla u(x, y) = f(x,y)$, in the
domain $[0,1]\times [0,1]$ with Dirichlet boundary conditions.
We choose the tensor coefficient $\boldsymbol{K}(x, y)=QMQ^{T}$, where $Q=\begin{bmatrix}
\cos(\theta) & -\sin(\theta)
\\ \sin(\theta) & \cos(\theta)
 \end{bmatrix}$, and $M=\begin{bmatrix}
\delta & 0
\\ 0 & 1
 \end{bmatrix}$. The parameters $0<\delta \leq 1$ and $\theta$ specify the strength
 and direction of anisotropy in the problem, respectively, with $\theta = 0$ giving the anisotropic problem $-\delta u_{xx} - u_{yy} = f$ and $\theta = \pi/2$ giving  $- u_{xx} - \delta u_{yy} = f$.  For $0 < \theta < \pi/2$, the axis of the small diffusion coefficient in the problem rotates clockwise from being in the positive $x$-direction for small $\theta$ to the positive $y$-direction for $\theta \approx \pi/2$.  Anisotropic problems cause difficulty for the greedy coarsening algorithm~\cite{maclachlan2007greedy}, where large grid and operator complexities and poor algorithmic performance were overcome by augmenting the coarse grids with the second pass of the Ruge-St\"{u}ben coarsening algorithm and using classical AMG interpolation in place of the AMGr interpolation operator. We emphasize that both the FD and FE discretizations of this problem result in non-diagonally dominant system matrices and, consequently, the convergence guarantee from~\cref{thm:amgr-cgsa} does not apply in this setting.  We include these problems, nonetheless, both to stress-test our approach and highlight the need for further research in this direction.

\begin{figure}[!ht]
\centering
\includegraphics[width=0.45\textwidth]{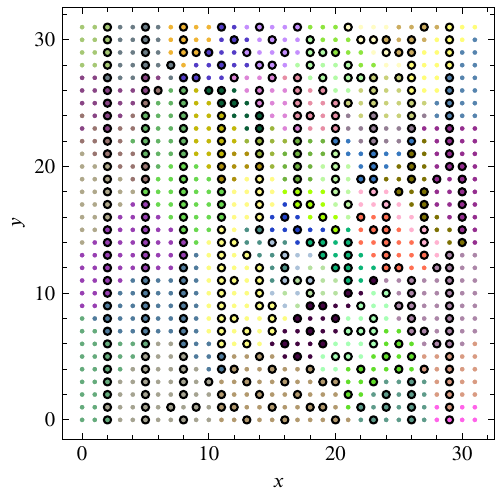} % chktex 8 chktex 29
\qquad
\includegraphics[width=0.45\textwidth]{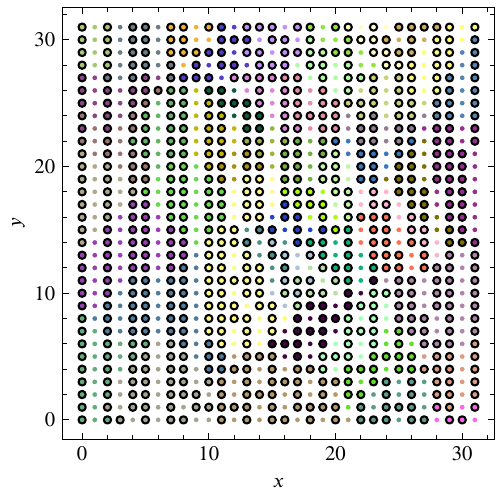} % chktex 8 chktex 29
\caption{Grid partitioning for an anisotropic diffusion problem with $\delta=10^{-6}$ and $\theta=\pi/3$, discretized on a $32\times 32$ mesh using the bilinear FE stencil.  The initial partitioning, at left, is generated using algebraic subdomains averaging 20 points per subdomain, using \num{2000000} SA steps per DoF and 1 SA step per DoF per GS sweep, and has 664 $F$-points and 360 $C$-points.  At right is the grid augmented using the second pass of the classical AMG algorithm, resulting in 343 $F$-points and 681 $C$-points.}\label{fig:cfnodes_struc_ani_32X32_FE_lloyd_20-51_2M_1_ba_2levpass}
\end{figure}

As an example of a problem in this class, we fix $\delta = 10^{-6}$ and $\theta = \pi/3$.  \Cref{fig:cfnodes_struc_ani_32X32_FE_lloyd_20-51_2M_1_ba_2levpass} shows the coarse-grid points selected for the bilinear finite-element discretization of the problem on a uniform $32\times 32$ mesh.  At left, we see that the partitioning generated by the SA algorithm correctly detects that this is an anisotropic problem with strong coupling primarily in the $x$-direction and weak coupling primarily in the $y$-direction, producing grids that are consistent with semi-coarsening, with some regions of coarsening along diagonals.  Unfortunately, however, this grid leads to poor convergence using AMGr interpolation.  As in the original experiments by MacLachlan and Saad~\cite{maclachlan2007greedy}, each fine-grid point has multiple coarse-grid neighbours (leading to an increase in operator complexity, due to the nonzero structure of $D_{\ssff}^{-1}A_{\ssfc}$), but since $A$ is no longer diagonally dominant in the anisotropic case, the assumption that $\begin{bmatrix} D_{\ssff} & -A_{\ssfc} \\ -A_{\ssfc}^T & A_{\sscc}\end{bmatrix}$ is positive semi-definite is violated, and the resulting performance is poor.  To overcome this failure, we augment the coarse-grid set using the second-pass algorithm from Ruge-St\"{u}ben AMG, using the classical strength of connection parameter of $0.30$ to determine strong connections in the graph.  At right of~\cref{fig:cfnodes_struc_ani_32X32_FE_lloyd_20-51_2M_1_ba_2levpass}, we show the resulting graph, which now satisfies the requirement that every pair of strongly connected $F$-points has a common $C$-neighbour (which is clearly not satisfied in the partitioning at left).  Unfortunately, this comes at a heavy cost, as the grid at right now has many more $C$-points than we would like.  Below, we verify that these grids lead to effective multigrid hierarchies, when coupled with classical AMG interpolation, but note the increased grid and operator complexities in these hierarchies.  A key question for future work (addressed in~\cite{TZaman_etal_2022a}) is whether we can determine algebraic interpolation operators for grids such as those at left that lead to effective AMGr performance.

\Cref{tab:conv_anisot_struc_geom-lloyd_1e-6_pi3_2lev-cgsa,tab:conv_anisot_struc_lloyd_1e-6_pi3_3lev_V_W-cgsa} present convergence factors and grid and operator complexities for two- and three-level cycles, respectively.  For geometric partitioning, we use \num{200000} SA steps per DoF with one or two SA steps per DoF per cycle and $5\times 5$ or $6\times 6$ subdomains (depending on what worked best for a given problem, reported in~\cref{tab:table5} in~\cref{appendix:data}).  Similar choices are made using algebraic partitioning, although we see slight improvements using \num{2000000} SA steps per DoF for some problems.  For the three-level tests, we uniformly use \num{200000} SA steps per DoF, with 1 SA step per DoF per GS sweep, and algebraic subdomains with average size of 36 points.  In~\cref{tab:conv_anisot_struc_geom-lloyd_1e-6_pi3_2lev-cgsa}, we again see that there is relatively little difference in results generated using geometric and algebraic choices of the Gauss-Seidel subdomains.  Notably, both choices lead to effective cycles for both the finite-difference and finite-element discretizations, albeit at the cost of increased grid and operator complexities.  \Cref{tab:conv_anisot_struc_lloyd_1e-6_pi3_3lev_V_W-cgsa} again shows degradation in performance moving from two-grid to three-grid V-cycles, but that three-grid W-cycles mostly recover comparable convergence to the two-level case, with some notable degradation for the finite-element operator.
\begin{table}[!ht]
	\caption{Performance of two-level AMGr for the anisotropic diffusion problem with $\delta = 10^{-6}$ and $\theta = \pi/3$ on structured meshes.}
	\centering
	%\begin{adjustbox}{width=0.6\textwidth}
	    %\resizebox{\textwidth}{!}{\begin{tabular}{|c|*{8}{c|}}
	    %\small\begin{tabular*}{\textwidth}{l@{\extracolsep{\fill}}*{9}{c}}
		{\begin{tabular}{cc|ccc|ccc}
			\toprule
			   & & \multicolumn{3}{c|}{Geometric Partitioning} & \multicolumn{3}{c}{Algebraic Partitioning} \\
      Scheme & Grid size& $\rho$ & $C_{\text{grid}}$ & $C_{\text{op}}$ & $\rho$ & $C_{\text{grid}}$ & $C_{\text{op}}$ \\
			\midrule
			   & $32\times32$ & 0.58 & 1.70 & 2.09 & 0.57 & 1.70 & 2.09 \\
			FD & $64\times64$ & 0.58 & 1.71 & 2.11 & 0.58 & 1.71 & 2.11 \\
			   & $128\times128$ & 0.60 & 1.72 & 2.12 & 0.63 & 1.72 & 2.12 \\
			   %& $256\times256$ &  &  &  &  &  & \\
			\midrule
			   & $32\times32$ & 0.62 & 1.65 & 2.02 & 0.61 & 1.67 & 2.00 \\
			FE & $64\times64$ & 0.61 & 1.67 & 2.09 & 0.62 & 1.67 & 2.08 \\
			   & $128\times128$ & 0.61 & 1.67 & 2.11 & 0.61 & 1.67 & 2.12 \\
			   %& $256\times256$ &  &  &  &  &  & \\
			\bottomrule
		\end{tabular}}\label{tab:conv_anisot_struc_geom-lloyd_1e-6_pi3_2lev-cgsa}
\end{table}
\begin{table}[!ht]
	\caption{Performance of three-level AMGr for the anisotropic diffusion problem with $\delta = 10^{-6}$ and $\theta = \pi/3$ on structured meshes.}
	\centering
	%\begin{adjustbox}{width=0.6\textwidth}
	    %\resizebox{\textwidth}{!}{\begin{tabular}{|c|*{8}{c|}}
	    %\small\begin{tabular*}{\textwidth}{l@{\extracolsep{\fill}}*{9}{c}}
		{\begin{tabular}{cc|cccc}
			\toprule
      Scheme & Grid size  & $\rho_V$ & $\rho_W$ & $C_{\text{grid}}$ & $C_{\text{op}}$ \\
			\midrule
			   & $32\times32$ & 0.72 & 0.57 & 2.16 & 3.13 \\
			FD & $64\times64$ & 0.76 & 0.60 & 2.16 & 3.20 \\
			   & $128\times128$ & 0.80 & 0.66 & 2.18 & 3.26 \\
			   %& $256\times256$ &  &  &  &  \\
			\midrule
			   & $32\times32$ & 0.74 & 0.61 & 2.06 & 2.86 \\
			FE & $64\times64$ & 0.81 & 0.69 & 2.07 & 3.05 \\
			   & $128\times128$ & 0.89 & 0.82 & 2.05 & 3.10 \\
			   %& $256\times256$ &  &  &  & \\
			\bottomrule
		\end{tabular}}\label{tab:conv_anisot_struc_lloyd_1e-6_pi3_3lev_V_W-cgsa}
\end{table}

%%%%%%%%%%%%%%%%%%%%%%%%%%%%%%%%%%%%%%%%%%%%%%%%%%%%%%%%%%%%%%%%%%%%%%%%%%%%%%

\subsection{Discretizations on unstructured grids}\label{subsec:unstructured-cgsa}

\begin{figure}[!ht]
\centering
\includegraphics{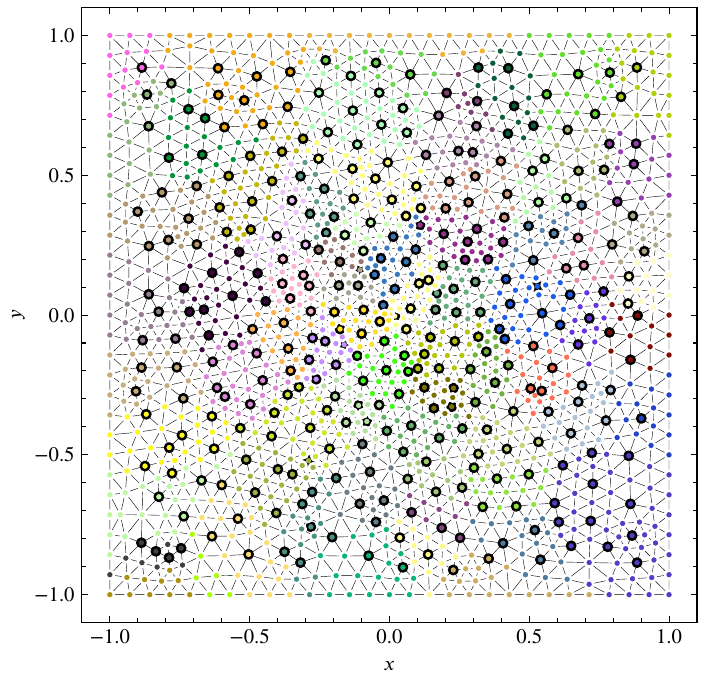}
\caption{Partitioning for the isotropic problem on the unstructured mesh with \num{1433} points, generated using \num{1000000} SA steps per DoF, 5 SA steps per DoF per GS sweep, and subdomains with an average of 20 points per subdomain.}\label{fig:cfnodes_unstruc_1433_FE_lloyd_20_71_1000000_5}
\end{figure}

We next consider two diffusion problems with homogeneous Dirichlet boundary conditions on unstructured triangulations of square domains, discretized using linear finite elements.  First, we consider an isotropic diffusion operator ($-\nabla\cdot\nabla u = f$) on $[-1,1]^2$, generated by taking an initially unstructured mesh, refining it a set number of times and, then, smoothing the resulting mesh.  We consider three levels of refinement for this example, resulting in meshes with \num{1433}, \num{5617}, and \num{22241} DoFs.  The SA-based partitioning scheme appears to perform similarly well in this setting, with a splitting found for the smallest mesh shown in~\cref{fig:cfnodes_unstruc_1433_FE_lloyd_20_71_1000000_5}.  While we no longer have hand-generated estimates of the optimal coarsening factors, we note that the SA-based partitioning scheme outperforms the greedy algorithm, generating $F$-sets with \num{1106}, \num{4241}, and \num{16604} points, for these three grids, respectively, in comparison to $F$-sets of \num{1024}, \num{3916}, and \num{15079} points. As above, this improvement comes at a cost: the offline time to partition the \num{5617} DoF problem using \num{200000} SA steps per DoF is over 12 hours for the two-level scheme.  Two- and three-level convergence factors, as well as grid and operator complexities for these systems are given in~\cref{tab:conv_isot_unstruc_lloyd_2-3lev_W-cgsa}.  Here, we observe that these results are quite similar to those seen for the structured-grid discretizations above, with some degradation in convergence seen for the three-level V-cycle results, but not in the W-cycle results.   Here, the three-level results are computed using \num{200000} SA steps per DoF, with 1 SA step per DoF per GS sweep, on subdomains with an average of 36 points per subdomain.  Slight modifications to these parameters yield small improvements in two-level results; as a result, we use these choices in the table and report them in~\cref{tab:tables78} in~\cref{appendix:data}.
%
%
% No. of F-points from greedy and SA are tabulated below.
% Unstructured FE (isotropic)
% #DoFs   F-points from greedy alg.     F-points from SA (max we got)
%----------------------------------------------------------------------
% 1433           1024                          1106
% 5617           3916                          4241
% 22241          15079                         16604

% Dom2
% #DoFs   F-points from greedy alg.     F-points from SA (max we got)
%---------------------------------------------------------------------
% 798            480                           524
% 3109           1722                          1947
% 12273          6527                          7391
%
%
%
%
% 2-level AMGr
%DoF                      2-level V
%-------------------------------------------------------------------------
% 1433       SAsteps/DoF:1M,SAsteps/DoF/cyc:5,subdom:20npa
% 5617       SAsteps/DoF:2M,SAsteps/DoF/cyc:2,subdom:36npa
% 22241      SAsteps/DoF:200K,SAsteps/DoF/cyc:1,subdom:36npa

% 3-level AMGr
%DoF                          3-level V                                          3-level W
%--------------------------------------------------------------------------------------------------------------------------
% 1433       SAsteps/DoF:200K,SAsteps/DoF/cyc:1,subdom:36npa     SAsteps/DoF:200K,SAsteps/DoF/cyc:1,subdom:36npa
% 5617       SAsteps/DoF:200K,SAsteps/DoF/cyc:1,subdom:36npa     SAsteps/DoF:200K,SAsteps/DoF/cyc:1,subdom:36npa
% 22241      SAsteps/DoF:200K,SAsteps/DoF/cyc:1,subdom:36npa     SAsteps/DoF:200K,SAsteps/DoF/cyc:1,subdom:36npa
%
%
\begin{table}[!ht]
	\caption{Performance of two- and three-level AMGr for isotropic problem on unstructured meshes.}
	\centering
	%\begin{adjustbox}{width=0.6\textwidth}
	    %\resizebox{\textwidth}{!}{\begin{tabular}{|c|*{10}{c|}}
	    %\resizebox{\textwidth}{!}{\begin{tabular}{|c||c|c|c||c|c|c||c|c|c|}
	    %\small\begin{tabular*}{\textwidth}{l@{\extracolsep{\fill}}*{9}{c}}
		{\begin{tabular}{c|ccc|cccc}
			\toprule
			      & \multicolumn{3}{c|}{Two-level cycle} & \multicolumn{4}{c}{Three-level cycles}\\
      \#DoF & $\rho$ & $C_{\text{grid}}$ & $C_{\text{op}}$ & $\rho_V$ & $\rho_W$ & $C_{\text{grid}}$ & $C_{\text{op}}$ \\
			\midrule
			            1433 & 0.66 & 1.23 & 1.31 & 0.78 & 0.65 & 1.28 & 1.39 \\
			            5617 & 0.71 & 1.25 & 1.32 & 0.79 & 0.70 & 1.30 & 1.42 \\
			           22241 & 0.75 & 1.25 & 1.33 & 0.84 & 0.75 & 1.32 & 1.45 \\
			\bottomrule
		\end{tabular}}\label{tab:conv_isot_unstruc_lloyd_2-3lev_W-cgsa}
\end{table}

\begin{figure}[!ht]
\centering
\includegraphics[width=0.45\textwidth]{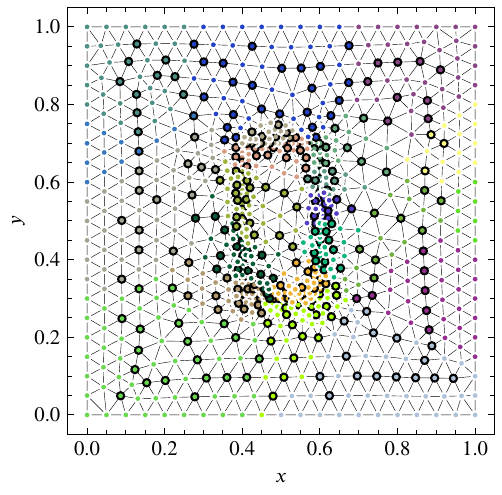}
\qquad
\includegraphics[width=0.45\textwidth]{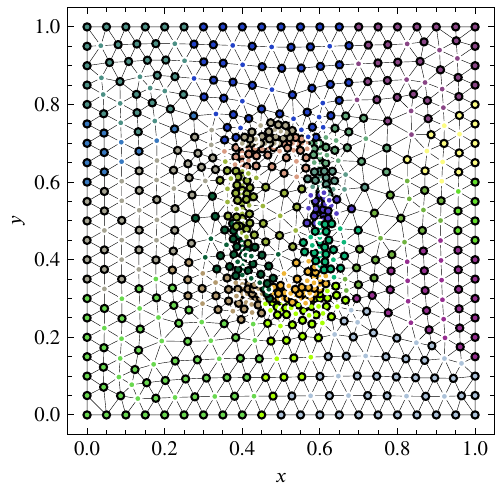}
\caption{Partitioning for the anisotropic problem on an unstructured mesh containing 798 points.  The partitioning at left was generated using \num{500000} SA steps per DoF, with 1 SA step per DoF per GS sweep, on subdomains with an average size of 36 points per subdomain, yielding 524 $F$-points and 274 $C$-points.  At right, this partitioning is augmented by the second pass of classical AMG coarsening, resulting in 263 $F$-points and 535 $C$-points.}\label{fig:cfnodes_unstruc_dom2_0p01_pi6_798_FE_lloyd_36_22_500000_1_ba_2levpass}
\end{figure}

We next consider the anisotropic diffusion operator (again with Dirichlet boundary conditions) considered above, with $\delta = 0.01$ and $\theta = \pi/3$, on an unstructured triangulation of the unit square taken from Brannick and Falgout~\cite{JJBrannick_RDFalgout_2010a}, matching the problem labelled 2D-M2-RLap in that paper.  As in Brannick and Falgout~\cite{JJBrannick_RDFalgout_2010a}, we consider three levels of refinement of the mesh, with \num{798}, \num{3109}, and \num{12273} DoFs, respectively. \Cref{fig:cfnodes_unstruc_dom2_0p01_pi6_798_FE_lloyd_36_22_500000_1_ba_2levpass} shows two different partitions for this mesh.  At left, we give the partitioning generated using the SA-based partitioning algorithm proposed here, and at right we give the splitting after a second pass of the Ruge-St\"{u}ben coarsening algorithm where strength of connection is computed using the classical strength parameter of $0.55$.  As with the structured-grid case above, we found the second pass is needed to achieve good convergence factors for the anisotropic problem, but that it does so at the expense of coarsening at a much slower rate.  Convergence factors, grid complexities, and operator complexities for these problems are reported in~\cref{tab:conv_anisot_unstruc_lloyd_2-3lev_W-cgsa}.  As above, the complexities are higher for the anisotropic operators than the isotropic (due to the use of the second pass).  Here, we observe  degradation in convergence with grid refinement, although again with less degradation for W-cycles than V-cycles.  The three-level results are computed with the same parameters as for the isotropic problem above, with more SA steps per DoF used for the two-level results, as this yielded slight improvements.  While the performance reported in~\cref{tab:conv_anisot_unstruc_lloyd_2-3lev_W-cgsa} is far from optimal, we note that the convergence factors for the larger two meshes are better than those reported for compatible relaxation~\cite{JJBrannick_RDFalgout_2010a}, while the operator complexities are comparable to those reported therein for BoomerAMG~\cite{VEHenson_UMYang_2002a}. We again emphasize that these anisotropic diffusion problems do not satisfy the convergence bounds stated above and, consequently, are included to stress-test the framework presented here; improved results for this problem are presented in~\cite{TZaman_etal_2022a}.
% Note regarding definition of angle:
% The matrix was formed using the code
% A0, b = fem_ani.gradgradform(mesh, kappa=kappa, f=f, degree=1)
% A0, b = fem.applybc(A0, b, mesh, bc), where, kappa was defined as:
% def kappa(x, y):
    %a = np.cos(theta)**2 + epsilon*np.sin(theta)**2
    %b = epsilon*np.cos(theta)**2 + np.sin(theta)**2
    %c = (1-epsilon)*np.cos(theta)*np.sin(theta)
    %Q = np.array([[a,  c],
                  %[c,  b]])
    %return Q
 % here, theta = sp.pi/6 was used.
% This gives the same matrix as the matrix that Dr. Rob gave us. However, in their paper, they mentioned that they used \theta=pi/3.
%
% Since we have switched orientations from M_{22} = \delta to M_{11} = \delta, I've switch the text to state \theta = \pi/3
%
%
%
% 2-level AMGr
%DoF                      2-level V
%-------------------------------------------------------------------------
% 798        SAsteps/DoF:500K,SAsteps/DoF/cyc:1,subdom:36npa
% 3109       SAsteps/DoF:2M,SAsteps/DoF/cyc:2,subdom:36npa
% 12273      SAsteps/DoF:200K,SAsteps/DoF/cyc:1,subdom:36npa

% 3-level AMGr
%DoF                          3-level V                                          3-level W
%--------------------------------------------------------------------------------------------------------------------------
% 798        SAsteps/DoF:200K,SAsteps/DoF/cyc:1,subdom:36npa     SAsteps/DoF:200K,SAsteps/DoF/cyc:1,subdom:36npa
% 3109       SAsteps/DoF:200K,SAsteps/DoF/cyc:1,subdom:36npa     SAsteps/DoF:200K,SAsteps/DoF/cyc:1,subdom:36npa
% 12273      SAsteps/DoF:200K,SAsteps/DoF/cyc:1,subdom:36npa     SAsteps/DoF:200K,SAsteps/DoF/cyc:1,subdom:36npa
%
%
\begin{table}[!ht]
	\caption{Performance of two-level and three-level AMGr for anisotropic problem on unstructured meshes.}
	\centering
	%\begin{adjustbox}{width=0.6\textwidth}
	    %\resizebox{\textwidth}{!}{\begin{tabular}{|c|*{8}{c|}}
	    %\resizebox{\textwidth}{!}{\begin{tabular}{|c||c|c|c||c|c|c||c|c|c|}
	    %\small\begin{tabular*}{\textwidth}{l@{\extracolsep{\fill}}*{9}{c}}
		{\begin{tabular}{c|ccc|cccc}
			\toprule
			      & \multicolumn{3}{c|}{Two-level cycle} & \multicolumn{4}{c}{Three-level cycles} \\
     \#DoF  & $\rho$ & $C_{\text{grid}}$ & $C_{\text{op}}$ & $\rho_V$ & $\rho_W$ & $C_{\text{grid}}$ & $C_{\text{op}}$ \\
     \midrule
			            %798 & 0.992 & 1.343 & 1.719 & 0.992 & 1.495 & 2.494 & 0.992 & 1.495 & 2.494 \\
			            %798 & 0.884 & 1.603 & 1.864 & 0.883 & 1.996 & 2.567 & 0.884 & 1.996 & 2.567 \\
			             798 & 0.69 & 1.67 & 1.82 & 0.84 & 0.71 & 2.11 & 2.41 \\
			           %3109 & 0.996 & 1.374 & 1.758 & 0.996 & 1.554 & 2.756 & 0.996 & 1.554 & 2.756 \\
			           %3109 & 0.921 & 1.603 & 1.906 & 0.918 & 1.996 & 2.694 & 0.920 & 1.996 & 2.694 \\
			            3109 & 0.75 & 1.67 & 1.87 & 0.84 & 0.75 & 2.09 & 2.52\\
			          %12273 & 0.997 & 1.398 & 1.807 & 0.997 & 1.600 & 2.950 & 0.997 & 1.600 & 2.950 \\
			          %12273 & 0.948 & 1.590 & 1.920 & 0.948 & 1.970 & 2.722 & 0.948 & 1.970 & 2.722 \\
			           12273 & 0.82 & 1.67 & 1.89 & 0.90 & 0.83 & 2.08 & 2.56 \\
			\bottomrule
		\end{tabular}}\label{tab:conv_anisot_unstruc_lloyd_2-3lev_W-cgsa}
\end{table}

%%%%%%%%%%%%%%%%%%%%%%%%%%%%%%%%%%%%%%%%%%%%%%%%%%%%%%%%%%%%%%%%%%%%%%%%%%%%%%

\subsection{Convection-diffusion problems on structured grids}

\begin{figure}[!ht]
\centering
\includegraphics[width=0.45\textwidth]{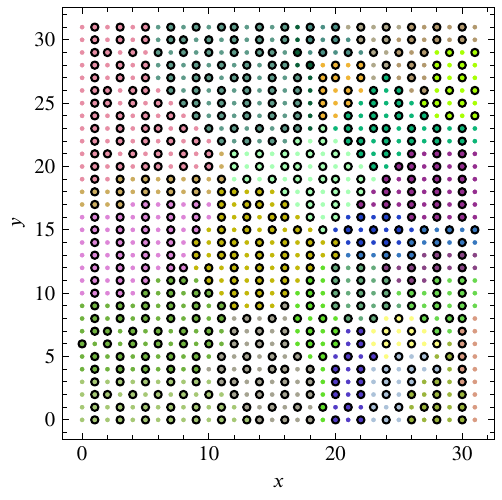}
%\hfill
\qquad
\includegraphics[width=0.45\textwidth]{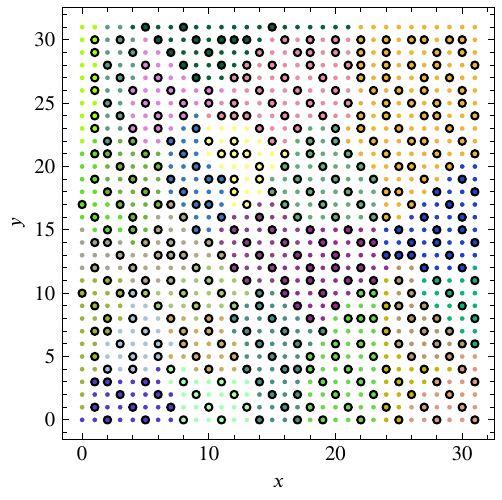}
\caption{Partitioning for convection-diffusion problems on uniform grids, for the grid-aligned case, at left, and the non-grid-aligned case, at right.  In both cases, we depict the partitioning for $\varepsilon=10^{-5}$ and $N=32$, generated using \num{200000} SA steps per DoF.}\label{fig:convection_diffusion}
\end{figure}

For our final tests, we consider the upwind finite-difference discretization of two singularly perturbed convection-diffusion equations.  While these problems lead to non-symmetric discretization matrices (and, as such, \Cref{thm:amgr-cgsa} no longer applies), they represent a plausible set of test problems for reduction-based methods due to their ``nearly triangular'' M-matrix structure~\cite{doi:10.1137/18M1193761}.  In particular, we consider the solution of the convection-diffusion equation,
\[
-\varepsilon\Delta u + \vec{b}\cdot\nabla u = f,
\]
with homogeneous Dirichlet boundary conditions, for two choices of the convection direction:
\begin{align*}
  \text{grid-aligned convection } & \vec{b} = \begin{bmatrix} 1 \\ 0 \end{bmatrix}, \\
  \text{non-grid-aligned convection } & \vec{b} = \begin{bmatrix} 2 \\ 3 \end{bmatrix}.
\end{align*}
We discretize these problems on uniform $N\times N$ meshes with the standard 5-point finite-difference for the diffusion term, and a first-order upwind finite-difference discretization for the convection terms.  In all experiments, we use \num{200000} SA steps per DoF to generate the partitioning, and generate the AMGr interpolation operator by computing $\begin{bmatrix} D_{\ssff}^{-1}A_{\ssfc} \\ I \end{bmatrix}$. We construct the restriction operator by taking $\begin{bmatrix} A_{\sscf}D_{\ssff}^{-1} & I \end{bmatrix}$.  \Cref{tab:convection_diffusion-cgsa} reports measured asymptotic convergence factors, along with grid and operator complexities for these problems, showing insensitivity to both problem size and the singular perturbation parameter, $\varepsilon$.  \Cref{fig:convection_diffusion} shows two sample partitioning, one for each convection direction, illustrated at $\varepsilon=10^{-5}$ for $N=32$.

\begin{table}[!ht]
	\caption{Performance of two-level AMGr for convection-diffusion problems on structured meshes.}
	\centering
		{\begin{tabular}{cc|ccc|ccc}
			\toprule
			      & & \multicolumn{3}{c|}{Grid-Aligned} & \multicolumn{3}{c}{Non-Grid-Aligned}\\
     $\varepsilon$ & $N$  & $\rho$ & $C_{\text{grid}}$ & $C_{\text{op}}$ & $\rho$ & $C_{\text{grid}}$ & $C_{\text{op}}$ \\
     \midrule
         & 16 & 0.617 & 1.50 & 1.93 & 0.617 & 1.34 & 1.65 \\
$10^{-3}$ & 32 & 0.617 & 1.50 & 2.01 & 0.617 & 1.36 & 1.71 \\
         & 64 & 0.617 & 1.51 & 2.02 & 0.617 & 1.36 & 1.74 \\
     \midrule
         & 16 & 0.617 & 1.50 & 1.93 & 0.617 & 1.34 & 1.63 \\
$10^{-4}$ & 32 & 0.617 & 1.50 & 1.97 & 0.617 & 1.36 & 1.71 \\
         & 64 & 0.617 & 1.51 & 2.03 & 0.617 & 1.36 & 1.74 \\
     \midrule
         & 16 & 0.617 & 1.50 & 1.92 & 0.617 & 1.33 & 1.63 \\
$10^{-5}$ & 32 & 0.617 & 1.51 & 2.00 & 0.617 & 1.36 & 1.71 \\
         & 64 & 0.617 & 1.51 & 2.03 & 0.617 & 1.36 & 1.73 \\
     \bottomrule
		\end{tabular}}\label{tab:convection_diffusion-cgsa}
\end{table}

\section{Conclusions and future work}\label{sec:conclusion-cgsa}

While existing heuristics for coarse-grid selection within AMG offer
reliable performance for a wide class of problems, there remains much
interest in both improving our understanding of AMG convergence and
developing new approaches that offer guarantees of convergence for
even wider classes of problems.  We believe a promising, yet
under-explored, area of research is in the extension of reduction-based
AMG approaches, that have been shown to be effective for both
isotropic diffusion
problems~\cite{SPMacLachlan_TAManteuffel_SFMcCormick_2006a, maclachlan2007greedy} and more interesting classes of
problems, such as some hyperbolic PDEs~\cite{doi:10.1137/18M1193761}.
In this paper, we propose a new coarsening algorithm for AMGr, based on applying simulated annealing to the optimization problem first posed by MacLachlan and Saad~\cite{maclachlan2007greedy}.
We choose SA as it is a widely used optimization technique that
can be applied to combinatorial optimization problems.  Other
techniques are certainly possible (e.g., genetic algorithms or
particle swarm optimization), and it may be that
those techniques lead to improved performance.
For isotropic problems on both structured and unstructured meshes, we find that the SA-based partitioning approach outperforms the original greedy algorithm, sometimes dramatically, while sharing some of the existing limitations of the AMGr framework, particularly for anisotropic problems.
We believe that the performance difference between the greedy and SA-based algorithm
is, simply, due to the
complex nature of the optimization problem at hand.  The greedy
algorithm makes the best ``local'' choice at each stage, but those
local choices force poor global configurations.  In particular,
the greedy approach has no opportunity to undo a choice already
made, while the SA approach allows that to happen.
Nonetheless, we see this as an important proof-of-concept, showing that randomized search and other derivative-free optimization algorithms can be successfully applied to the combinatorial optimization problems in AMG coarsening.

The main drawback of this approach is the high computational cost of simulated annealing.  For example, computing a single coarse grid for a $32\times 32$ mesh with \num{200000} SA steps per DoF required around 2.5 hours on a modern workstation, primarily for evaluating the fitness functional in the optimization, with expected scaling in problem size and number of SA steps per DoF.  Two key questions that this raises are whether the fitness functional can be approximated in a more efficient manner (e.g., using a value neural network) and whether other optimization techniques that rely on fewer samplings of the fitness functional can be applied.  Both of these are topics of our current research.

This research also exposes a known weakness of the AMGr methodology (and, to our knowledge, one of all AMG-like algorithms with guarantees on convergence rates) when applied to problems that are not diagonally dominant (such as finite-element discretization of anisotropic diffusion equations).  Another current research direction is identifying whether the diagonal choice of $D_{\ssff}$ can be generalized to yield better convergence (while retaining a guaranteed convergence rate), and whether such changes can be accommodated within the optimization problems solved here. Preliminary results in this direction are presented in~\cite{TZaman_etal_2022a}.

%%%%%%%%%%%%%%%%%%%%%%%%%%%%%%%%%%%%%%%%%%%%%%%%%%%%%%%%%%%%%%%%%%%%%%%%%%%%%%

\section*{Acknowledgments}

The work of S.P.M. was partially supported by an NSERC Discovery Grant. This work does not have any conflicts of interest.

\bibliographystyle{elsarticle-num-names}
\bibliography{refs}

%%%%%%%%%%%%%%%%%%%%%%%%%%%%%%%%%%%%%%%%%%%%%%%%%%%%%%%%%%%%%%%%%%%%%%%%%%%%%%

\appendix
%\newpage

\section{Additional Data}\label{appendix:data}
In this section, we provide some additional tables to provide more complete data on the experiments reported above.
%%%%%%%%%%%%%%%%%%%%%%%%%%%

\begin{table}[!ht]
  \centering
  \begin{tabular}{c|ccccc}
    \toprule
    SA steps / DoF & $2\times 2$ & $3\times 3$ & $4\times 4$ & $5\times 5$ & $6\times 6$ \\
    \midrule
    \num{5000}    & 50 & 20 & 100 &  2 & 1 \\
    \num{200000}  & 10 &  1 &  25 &  5 & 2 \\
    \num{2000000} & -- & -- &  25 & 20 & 5 \\
    \bottomrule
\end{tabular}
  \caption{Number of SA steps per DoF per sweep used for numerical results in left-hand panel of~\cref{fig:maxfptVsubdomshape_o2L_maxfptVsastep_5spdpc_6X6_FD32_struc_geom}.}\label{tab:stepsFig3left}
\end{table}

%%%%%%%%%%%%%%%%%%%%%%%%%%%

\begin{table}[!ht]
  \centering
  \begin{tabular}{cccccccccc}
    \toprule
    SA steps / DoF\@: & \num{3000} & \num{5000} & \num{7500} & \num{10000} & \num{50000} & \num{100000} & \num{200000} & \num{500000} & \num{2000000}\\
    \midrule
    SA steps / DoF / cycle\@: & 1 & 1 & 1 & 2 & 1 & 2 & 2 & 2 & 5 \\
    \bottomrule
\end{tabular}
  \caption{Number of SA steps per DoF per sweep used for numerical results in right-hand panel of~\cref{fig:maxfptVsubdomshape_o2L_maxfptVsastep_5spdpc_6X6_FD32_struc_geom}.}\label{tab:stepsFig3right}
\end{table}

%%%%%%%%%%%%%%%%%%%%%%%%%%%
\begin{table}[!ht]
\small
  \centering
  \begin{tabular}{cc|cccccc}
  \toprule
Partitioning &  & \multicolumn{6}{c}{FD Discretization} \\
& & $8\times 8$ & $16\times 16$ & $32\times 32$ & $64\times 64$ & $128\times 128$ & $256\times 256$ \\
\midrule
     & SA steps per DoF & \num{1000} &  \num{10000} & \num{2000000} &  \num{2000000} &  \num{200000} & \num{200000} \\
   Geometric & SA steps per DoF per cycle & 1 & 100 & 5 & 1 & 1 & 1\\
    & Subdomain size & $2\times 2$ & $4\times 4$ & $6\times 6$ & $6\times 6$ & $6\times 6$ & $6\times 6$\\
    \hline
     & SA steps per DoF & -- & -- & \num{2000000} &  \num{2000000} &  \num{200000} & -- \\
   Algebraic & SA steps per DoF per cycle & -- & -- & 2 & 1 & 1 & -- \\
    & Subdomain size & -- & -- & 36 & 36 & 36 & -- \\
    \bottomrule
\end{tabular}
  \caption{Number of SA steps and SA steps per DoF per sweep used for numerical results in left-hand panel of~\cref{fig:maxfptVmeshsize_FDE_struc_geom_lloyd3}.} \label{tab:stepsFig7left}
\end{table}
%%%%%%%%%%%%%%%%%%%%%%%%%%%

\begin{table}[!ht]
\small
  \centering
  \begin{tabular}{cc|cccccc}
  \toprule
Partitioning &  & \multicolumn{6}{c}{FE Discretization} \\
& & $8\times 8$ & $16\times 16$ & $32\times 32$ & $64\times 64$ & $128\times 128$ & $256\times 256$ \\
\midrule
     & SA steps per DoF & \num{1000} &  \num{200000} & \num{2000000} &  \num{1000000} &  \num{200000} & \num{200000} \\
   Geometric & SA steps per DoF per cycle & 25 & 25 & 50 & 100 & 1 & 1\\
    & Subdomain size & $2\times 2$ & $3\times 3$ & $4\times 4$ & $3\times 3$ & $5\times 5$ & $5\times 5$\\
    \hline
     & SA steps per DoF & -- & -- & \num{2000000} &  \num{200000} &  \num{2000000} & -- \\
   Algebraic & SA steps per DoF per cycle & -- & -- & 1 & 1 & 1 & -- \\
    & Subdomain size & -- & -- & 25 & 25 & 25 & -- \\
    \bottomrule
\end{tabular}
  \caption{Number of SA steps and SA steps per DoF per sweep used for numerical results in right-hand panel of~\cref{fig:maxfptVmeshsize_FDE_struc_geom_lloyd3}.} \label{tab:stepsFig7right}
\end{table}

%%%%%%%%%%%%%%%%%%%%%%%%%%%

\begin{table}[!ht]
  \centering
  \begin{tabular}{c|ccccc}
    \toprule
    SA steps / DoF & $2\times 2$ & $3\times 3$ & $4\times 4$ & $5\times 5$ & $6\times 6$ \\
    \midrule
    \num{5000}    & 20 & 10 &   5 &  1 & 1 \\
    \num{200000}  & 10 &  2 &  50 & 20 & 2 \\
    \num{2000000} & -- & -- &  50 &  1 & 1 \\
    \bottomrule
\end{tabular}
  \caption{Number of SA steps per DoF per sweep used for numerical results in right-hand panel of~\cref{fig:maxfptVsastepPsubdomPcyc_FE32_5_maxfptVsubdomshape_struc_geom}.}\label{tab:stepsFig8}
\end{table}

%%%%%%%%%%%%%%%%%%%%%%%%%%%

\begin{table}[!ht]
\small
  \centering
  \begin{tabular}{c|cccccccc}
    \toprule
$\#$ of nodes per aggregate & 2 & 5 & 8 & 11 & 14 & 16 & 25 & 36\\
  \midrule
  FD\@: & & & & & & & & \\
  SA steps per DoF & \num{200000} & \num{7500} & \num{100000} &  \num{200000} & \num{200000} & \num{1000000} & \num{500000} & \num{2000000} \\
  SA steps per DoF per cycle & 10 & 10 & 20 & 5 & 50 & 25 & 1 & 2 \\
    \midrule
   FE\@: & & & & & & & & \\
  SA steps per DoF & \num{7500} & \num{200000} & \num{2000000} &  \num{200000} & \num{200000} & \num{200000} & \num{2000000} & \num{2000000} \\
  SA steps per DoF per cycle & 5 & 100 & 20 & 50 & 1 & 1 & 1 & 1 \\
  \midrule
\end{tabular}
\caption{Number of SA steps and SA steps per DoF per sweep used for numerical results in~\cref{fig:maxfptVsubdomshape-nnpagg_FDE32_struc_lloyd}.}\label{tab:stepsFig10}
\end{table}

%%%%%%%%%%%%%%%%%%%%%%%%%%%

\begin{table}[!ht]
  \centering
  \begin{tabular}{cc|ccc|ccc}
    \toprule
Partitioning & & \multicolumn{3}{c|}{FD Discretization} & \multicolumn{3}{c}{FE Discretization}\\
 & & $32\times 32$ & $64\times 64$ & $128\times 128$ &  $32\times 32$ & $64\times 64$ & $128\times 128$ \\
\midrule
    & SA steps per DoF & \num{200000} &  \num{200000} &  \num{200000} &  \num{200000} &  \num{200000} &  \num{200000}\\
    Geometric & SA steps per DoF per cycle & 1 & 1 & 1 & 2 & 1 & 1 \\
    & Subdomain size & $6\times 6$ &  $6\times 6$ &  $6\times 6$ &  $5\times 5$ &   $5\times 5$ &   $5\times 5$ \\
    \midrule
    & SA steps per DoF & \num{200000} &  \num{200000} &  \num{200000} &  \num{2000000} &  \num{200000} &  \num{2000000}\\
    Algebraic &  SA steps per DoF per cycle & 1 & 1 & 1 & 1 & 1 & 2 \\
    & Subdomain size & 16 & 36 & 36 & 20 & 36 & 36 \\
    \bottomrule
\end{tabular}
  \caption{Number of SA steps per DoF and SA steps per DoF per sweep and subdomain sizes used for two-level numerical results in~\cref{tab:conv_anisot_struc_geom-lloyd_1e-6_pi3_2lev-cgsa}.}\label{tab:table5}
\end{table}

%%%%%%%%%%%%%%%%%%%%%%%%%%%

\begin{table}[!ht]
  \centering
  \begin{tabular}{c|ccc|ccc}
    \toprule
 & \multicolumn{3}{c|}{\Cref{tab:conv_isot_unstruc_lloyd_2-3lev_W-cgsa}} & \multicolumn{3}{c}{\Cref{tab:conv_anisot_unstruc_lloyd_2-3lev_W-cgsa}}\\
 & 1433 & 5617 & 22241 & 798 & 3109 & 12273 \\
\midrule
    SA steps per DoF & \num{1000000} &  \num{2000000} &  \num{200000} &  \num{500000} &  \num{2000000} &  \num{200000}\\
    SA steps per DoF per cycle & 5 & 2 & 1 & 1 & 2 & 1 \\
    Subdomain size & 20 & 36 & 36 & 36 & 36 & 36 \\
    \bottomrule
\end{tabular}
  \caption{Number of SA steps per DoF and SA steps per DoF per sweep and subdomain sizes used for two-level numerical results in~\cref{tab:conv_isot_unstruc_lloyd_2-3lev_W-cgsa,tab:conv_anisot_unstruc_lloyd_2-3lev_W-cgsa}.}\label{tab:tables78}
\end{table}

%%%%%%%%%%%%%%%%%%%%%%%%%%%

\end{document}